\documentclass[11pt,oneside,english]{amsart}

\usepackage[T1]{fontenc}
\usepackage[latin9]{inputenc}
\usepackage[letterpaper]{geometry}
\usepackage{color}
\usepackage{float}
\usepackage{textcomp}
\usepackage{amsthm}
\usepackage{amssymb}
\usepackage{graphicx}
\usepackage{setspace}
\usepackage{esint}
\usepackage{nameref}

\makeatletter

\floatstyle{ruled}
\newfloat{algorithm}{tbp}{loa}
\providecommand{\algorithmname}{Algorithm}
\floatname{algorithm}{\protect\algorithmname}

\numberwithin{equation}{section}

\usepackage{ifpdf} 
\ifpdf 
 \IfFileExists{lmodern.sty}{\usepackage{lmodern}}{}
\fi 
\usepackage{color}

\usepackage{fancyhdr}
\fancypagestyle{plain}{\fancyfoot{}\fancyhead[C]{\text{\thepage}}}
\usepackage{graphicx}

\usepackage{hyperref}
\makeatother

\usepackage{babel}
\begin{document}

\title{Intrinsic Frequency Analysis and Fast Algorithms}

\author{Peyman Tavallali$^{1,*}$, Hana Koorehdavoudi$^{2}$, Joanna Krupa$^{3}$}
\begin{abstract}
Intrinsic Frequency (IF) has recently been introduced as an ample
signal processing method for analyzing carotid and aortic pulse pressure
tracings. The IF method has also been introduced as an effective approach
for the analysis of cardiovascular system dynamics. The physiological
significance, convergence and accuracy of the IF algorithm has been
establis\textcolor{black}{hed in prior works. In this paper, we show
that the IF method could be derived by appropriate mathematical approximations
from the Navier-Stokes and elasticity equations. We further introduce
a fast algorithm for the IF method based on the mathematical analysis
of this method. }In particular, we demonstrate that the IF algorithm
can be made faster, by a factor or more than 100 times, using a proper
set of initial guesses based on the topology of the problem, fast
analytical solution at each point iteration, and substituting the
brute force algorithm with a pattern search method. Statistically,
we observe that the algorithm presented in this article complies well
with its brute-force counterpart. Furthermore, we will show that on
a real dataset, the fast IF method can draw correlations between the
extracted intrinsic frequency features and the infusion of certain
drugs. \textcolor{black}{In general, this paper aims at a mathematical
analysis of the IF method to show its possible origins and also to
present faster algorithms.}
\end{abstract}

\maketitle
$^{1}$Division of Engineering and Applied Sciences, California Institute
of Technology, 1200 East California Boulevard, MC 205-45, Pasadena,
CA 91125, USA

$^{2}$Aerospace and Mechanical Engineering, University of Southern
California, Los Angeles, CA 90089-1453, USA

$^{3}$Avicena LLC, 2400 N Lincoln Ave, Altadena, CA 91001, USA

$^{*}$Corresponding Author, email: ptavalla@caltech.edu, tavallali@gmail.com

\section{Introduction}

Cardiovascular diseases (CVDs) and stroke are major causes of death
in the United States. The total cost related to CVDs and stroke was
estimated to be more than \$316 billion in 2011-2012 \cite{lloyd2010heart,mozaffarian2016executive}.
Hence, clinical measurements of cardiovascular health indices are
of great importance. These methods and measurements are essential
tools for monitoring cardiovascular health due to their relative availability.
For example, Left Ventricular Ejection Fraction (LVEF) is a measure
of left ventricular contractility \cite{borlaug2009contractility}
and Carotid-Femoral Pulse Wave Velocity (cfPWV) is a measure of aortic
stiffness \cite{laurent2006expert}. 

However, current methods of measuring such indices are expensive,
sometimes invasive, prone to measurement errors, and not necessarily
easy to use. For example, 2D LVEF echocardiography is not accurate
compared to more expensive and laborious gold standard cardiac MRI
method \cite{dewey2006evaluation,greupner2012head,hoffmann2005assessment,hoffmann2014analysis}.
As another example, obtaining accurate cfPWV measurements often requires
certain medical devices and a well-trained staff within a clinical
setting \cite{rajzer2008comparison}. Consequently, continuous measurement
of these indices is not practical. These limitations emphasize the
need for new cardiovascular monitoring methods.

Intrinsic Frequency (IF) has been established as a new method of cardiovascular
monitoring through a novel signal processing methodology \cite{pahlevan2014intrinsic}.
The IF method needs only an uncalibrated pulse pressure \cite{tavallali2015convergence}
to extract pertinent information regarding the cardiovascular health
of an individual \cite{pahlevan2014intrinsic}. The IF method has
also been shown to be capable of non-invasively measuring LVEF by
means of an iPhone camera \cite{pahlevan2017noninvasive}. We believe
that methods like IF are of clinical and financial benefit in addressing
cardiovascular monitoring.

\textcolor{black}{In this paper, at first, we provide an overview
of the IF method. Next, we present an approximate derivation of the
IF model by combining Navier-Stokes equations and continuity with
elasticity equations. This helps to build a solid mathematical foundation
for the IF method and the analysis that follows. Later, we analyze
the IF algorithm in the space of feasible solutions, and based on
that, we introduce a new version of the IF algorithm which is faster
than the current brute-force IF method }\cite{tavallali2015convergence}\textcolor{black}{{}
while maintaining the same accuracy. We then perform a case study
on real pressure waveforms drawn from canine data using our new algorithm.
We will see that the fast IF algorithm is capable of capturing the
effects of different drug infusions on a canine subject. }

\section{Brief Overview of IF method }

\subsection{A History of Analyzing Cardiovascular Pulse Waveform}

Blood pressure was first measured by Hales in 1735 \cite{hall1987stephen}.
In his measurements, he found that blood pressure is not constant
in the arterial system. He related these variations to the elasticity
of the arteries \cite{hall1987stephen}. Currently, it is known that
the shape of the arterial pulse wave is intimately related to the
physiology and pathology of the whole arterial system \cite{alastruey2012arterial}.
There has been much research on analyzing the dynamics of blood pressure
and flow in arterial systems \cite{avolio2009arterial,denardo2010pulse,nichols2011mcdonald,nichols2005mcdonald,wang2008ventricular}.
Specifically, there are two main approaches to analyzing cardiovascular
pulse wave data. One approach is based on a systematic mathematical
framework for the cardiovascular system. The other is based on directly
analyzing the pulse pressure waveform using signal processing methods.

An example of the systematic framework can be seen with the set of
Windkessel models \cite{westerhof1969analog}. The formulation of
a minimal lumped model of the arterial system was first presented
by Westerhof et al. \cite{westerhof1969analog}. Based on a Windkessel
model, the arterial system dynamics have been modeled through a combination
of different elements such as resistance, compliance and impedance.
In this simplified model of the arterial system, the blood flow dynamics
is modeled by the interaction between the elements (assuming the blood
flow acts as the current in the system). Because of the type of modeling,
the wave transmission of the blood flow is neglected. As a result,
the Windkessel models is not able to represent the entire dynamics
of the blood flow in an arterial system accurately. 

On the other hand, there are various methods for direct analysis of
an arterial pulse waveform, in both time and frequency domains \cite{nichols2005mcdonald}.
For example, the impedance method, which is based on Fourier transform,
is a common method to analyze the pressure waveform in the frequency
domain \cite{avolio2009arterial}. As an example, Milnor has shown
that the pressure and flow waveforms can be a superposition of several
harmonics using the Fourier method \cite{milnor1989hemodynamics}.
Another method to investigate the pressure wave in the time domain
is the wave intensity analysis which is based on wavelet transform
\cite{de2009blood}. These methods do not necessarily convey a physical
understanding of the cardiovascular system.

The IF algorithm presented in \cite{tavallali2015convergence} is
analyzing a pulse waveform through a direct time-frequency signal
processing machinery setting, from a quantitative perspective. Although,
in previous work \cite{pahlevan2014intrinsic}, we tried to qualitatively
express a systems approach to the IF formulation, the quantitative
picture has not yet been clear. However, in this article, we show
this connection from a quantitative perspective. 

\subsection{IF Formulation}

\textcolor{black}{In the IF method,} the aortic pressure waveform
at time $t\in\left[0,T\right)$, for a cardiac period $T$, can be
represented as
\begin{equation}
\begin{array}{ccc}
S\left(a_{i},b_{i},\bar{p},\omega_{i};t\right) & = & \left(a_{1}\cos\omega_{1}t+b_{1}\sin\omega_{1}t+\bar{p}\right)\mathbf{1}_{\left[0,T_{0}\right)}\left(t\right)+\\
 &  & \left(a_{2}\cos\omega_{2}t+b_{2}\sin\omega_{2}t+\bar{p}\right)\mathbf{1}_{\left[T_{0},T\right)}\left(t\right),
\end{array}\label{eq: The Trend}
\end{equation}
with a continuity condition at $T_{0}$ and periodicity at $T$. In
this formulation, the \emph{indicator function} is defined as 
\[
\mathbf{1}_{\left[x,y\right)}\left(t\right)=\left\{ \begin{array}{cc}
1, & x\leq t<y,\\
0, & else.
\end{array}\right.
\]
Here, $a_{1}$, $b_{1}$, $a_{2}$ and $b_{2}$ are the envelopes
of the IF model fit. $\omega_{1}$ and $\omega_{2}$ are the Intrinsic
Frequencies (IFs) of the waveform. Further, $\bar{p}$ is the mean
pressure during the cardiac cycle. \textcolor{black}{This type of
formulation embeds the coupling and decoupling of heart and aorta.}

The goal of the IF model (\ref{eq: The Trend}) is to extract a fit,
called Intrinsic Mode Function (IMF), that carries most of the energy
(information) from an observed pressure waveform $f\left(t\right)$
in one period. The latter is done by solving the following optimization
problem \cite{tavallali2015convergence}:
\begin{equation}
\begin{array}{cc}
\underset{a_{i},b_{i},\omega_{i},\bar{p}}{minimize} & \left\Vert f\left(t\right)-S\left(a_{i},b_{i},\bar{p},\omega_{i};t\right)\right\Vert _{2}^{2}\end{array}\label{eq: Continuous min}
\end{equation}
\begin{equation}
subject\,to\begin{array}{ccc}
a_{1}\cos\omega_{1}T_{0}+b_{1}\sin\omega_{1}T_{0} & = & a_{2}\cos\omega_{2}T_{0}+b_{2}\sin\omega_{2}T_{0},\\
a_{1} & = & a_{2}\cos\omega_{2}T+b_{2}\sin\omega_{2}T.
\end{array}\label{eq: Conditions}
\end{equation}
Here, $\left\Vert \right\Vert _{2}$ is the $L^{2}$-norm. The first
linear condition in this optimization enforces the continuity of the
extracted IMF at the dicrotic notch. The second one imposes the periodicity.
The mathematical convergence and accuracy of the IF algorithm have
been explained in a previous work \cite{tavallali2015convergence}.
In the next sections, we explore the foundation of the IF algorithm
and propose a faster IF algorithm.

\section{Approximate Derivation of the IF Model}

As mentioned earlier, in a previous work \cite{pahlevan2014intrinsic},
we tried to express a systems approach to the IF formulation qualitatively.
However, in this article, we show this connection from a quantitative
perspective. This section is devoted to this purpose. 

\begin{figure}[H]
\includegraphics[scale=0.4]{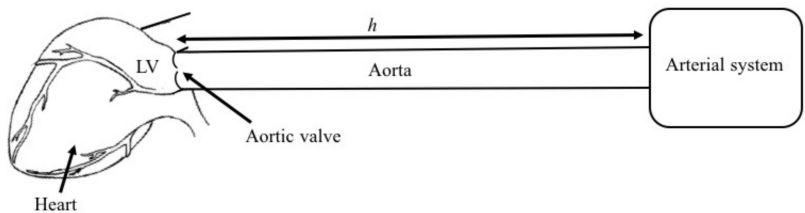}

\caption{Simplified cardiovascular system model schematic}

\label{Fig: Model}
\end{figure}

In this paper, we assume that the Left Ventricle (LV), the aortic
valve, aorta and the arterial system can be represented by a simplified
model as shown in Figure \ref{Fig: Model}. Here, the LV and the aortic
valve are assumed to be the boundary condition at the entrance of
the aortic tube and the arterial system is the terminal boundary condition
of the aortic tube. The boundary condition at the entrance of the
aortic tube changes from an LV boundary condition to a closed valve
boundary condition, at the dicrotic notch time $T_{0}$ during a cardiac
cycle $\left[0,T\right]$. \textcolor{black}{We further assume that
blood is a Newtonian incompressible fluid, the aorta is a straight
and sufficiently long elastic tube with a constant circular cross
section and there is no external force causing flow rotation. These
assumptions are not all satisfied in a real cardiovascular system.
However, they are useful in estimating the general behavior of blood
in aorta.}

Combining the Navier-Stokes equations and continuity with the elasticity
equation, \textcolor{black}{we can drive a model for the flow $Q\left(x,t\right)$
and the pressure $P\left(x,t\right)$ along the length $x$ of an
aorta as follow}

\textcolor{black}{
\begin{equation}
-\frac{\partial P}{\partial x}\left(x,t\right)=L\frac{\partial Q}{\partial t}\left(x,t\right)+RQ\left(x,t\right),\label{eq: Linearized Momentum}
\end{equation}
\begin{equation}
-\frac{\partial Q}{\partial x}\left(x,t\right)=C\frac{\partial P}{\partial t}\left(x,t\right).\label{eq: Simplified Elasticity}
\end{equation}
The step by step derivation of these equations is presented in \nameref{sec:Appendix-A}.
Parameters $L$, $R$, and $C$ represent inductance, resistance,
and compliance of the blood in aorta. Here, $0\leq x\leq h$, where
$h$ represents the aortic length. }This model has also been discussed
and simulated numerically in \cite{berger1994differential} with a
complex set of boundary conditions. Here, our main concentration will
be on the aortic tube oscillatory waveform solutions. Next, we will
show that we can derive (\ref{eq: The Trend}) from (\ref{eq: Linearized Momentum})
and (\ref{eq: Simplified Elasticity}).

Since the input to the IF model (\ref{eq: The Trend}) is a pressure
waveform, we need to extract an equation for the pressure \textcolor{black}{$P\left(x,t\right)$}
from Equations (\ref{eq: Linearized Momentum}) and (\ref{eq: Simplified Elasticity})
by eliminating the flow. Combining Equations (\ref{eq: Linearized Momentum})
and (\ref{eq: Simplified Elasticity}) results in 

\begin{equation}
CL\frac{\partial^{2}P}{\partial t^{2}}\left(x,t\right)+CR\frac{\partial P}{\partial t}\left(x,t\right)=\frac{\partial^{2}P}{\partial x^{2}}\left(x,t\right).\label{eq: Pressure Damped Wave}
\end{equation}
Taking $P\left(x,t\right)=\mathcal{K}\left(t\right)p\left(x,t\right)+\bar{p}$,
with $\bar{p}$ as the the mean pressure, we can write equation (\ref{eq: Pressure Damped Wave})
as

\begin{equation}
\begin{array}{c}
\left(CL\ddot{\mathcal{K}}\left(t\right)+CR\dot{\mathcal{K}}\left(t\right)\right)p\left(x,t\right)+\left(2CL\dot{\mathcal{K}}\left(t\right)+CR\mathcal{K}\left(t\right)\right)\frac{\partial p}{\partial t}\left(x,t\right)+CL\mathcal{K}\left(t\right)\frac{\partial^{2}p}{\partial t^{2}}\left(x,t\right)\\
=\mathcal{K}\left(t\right)\frac{\partial^{2}p}{\partial x^{2}}\left(x,t\right).
\end{array}\label{eq:44}
\end{equation}
Here, we have used the dot notation to represent the time derivative.
We can simplify the term in front of $\frac{\partial p}{\partial t}\left(x,t\right)$,
in (\ref{eq:44}), by setting $2CL\dot{\mathcal{K}}\left(t\right)+CR\mathcal{K}\left(t\right)=0$.
The latter has a solution $\mathcal{K}\left(t\right)=Ke^{-\frac{R}{2L}t}$
for some constant $K$. This reduces Equation (\ref{eq:44}) into
\begin{equation}
CL\frac{\partial^{2}p}{\partial t^{2}}\left(x,t\right)-\frac{CR^{2}}{4L}p\left(x,t\right)=\frac{\partial^{2}p}{\partial x^{2}}\left(x,t\right).\label{eq: KleinGordon}
\end{equation}
The solution of Equation (\ref{eq: KleinGordon}) can be expressed
in terms of eigenfunctions. In other words, using the method of separation
of the variables, one can express the solution of Equation (\ref{eq: KleinGordon})
as 
\begin{equation}
p\left(x,t\right)=\sum_{n=1}^{\infty}T_{n}\left(t\right)X_{n}\left(x\right)\label{eq:48}
\end{equation}
for

\begin{equation}
T_{n}\left(t\right)=\mathcal{\alpha}_{n}sin\left(\omega_{n}t\right)+\beta_{n}cos\left(\omega_{n}t\right),\label{eq:52}
\end{equation}
\begin{equation}
X_{n}\left(x\right)=\zeta_{n}sin\left(\sqrt{CL\left(\omega_{n}\right)^{2}-\frac{CR^{2}}{4L}}x\right)+\eta_{n}cos\left(\sqrt{CL\left(\omega_{n}\right)^{2}-\frac{CR^{2}}{4L}}x\right),\label{eq:53}
\end{equation}
and some constants $\mathcal{\alpha}_{n}$, $\beta_{n}$, $\zeta_{n}$
and $\eta_{n}$. As a result, the solution of (\ref{eq: Pressure Damped Wave})
can be expressed as 
\begin{equation}
\begin{array}{c}
P\left(x,t\right)=\bar{p}+Ke^{-\frac{R}{2L}t}\sum_{n=1}^{\infty}\Biggl\{\left(\mathcal{\alpha}_{n}sin\left(\omega_{n}t\right)+\beta_{n}cos\left(\omega_{n}t\right)\right)\\
\left(\zeta_{n}sin\left(\sqrt{CL\left(\omega_{n}\right)^{2}-\frac{CR^{2}}{4L}}x\right)+\eta_{n}cos\left(\sqrt{CL\left(\omega_{n}\right)^{2}-\frac{CR^{2}}{4L}}x\right)\right)\Biggr\}.
\end{array}\label{eq: Solution1}
\end{equation}
The variables $\omega_{n}$ can be expressed based on the boundary
conditions of the aortic tube. We need to emphasize that for a period
of the cardiac cycle $\left[0,T\right)$, the boundary conditions
change before and after the dicrotic notch $T_{0}$. Hence, for $t\in\left[0,T\right)$,
Equation (\ref{eq: Solution1}) can be written as
\begin{equation}
\begin{array}{c}
P\left(x,t\right)=\bar{p}\\
+\mathbf{1}_{\left[0,T_{0}\right)}\left(t\right)K^{1}e^{-\frac{R}{2L}t}\sum_{n=1}^{\infty}\Biggl\{\left(\mathcal{\alpha}_{n}^{1}sin\left(\omega_{n}^{1}t\right)+\beta_{n}^{1}cos\left(\omega_{n}^{1}t\right)\right)\\
\left(\zeta_{n}^{1}sin\left(\sqrt{CL\left(\omega_{n}^{1}\right)^{2}-\frac{CR^{2}}{4L}}x\right)+\eta_{n}^{1}cos\left(\sqrt{CL\left(\omega_{n}^{1}\right)^{2}-\frac{CR^{2}}{4L}}x\right)\right)\Biggr\}\\
+\mathbf{1}_{\left[T_{0},T\right)}\left(t\right)K^{2}e^{-\frac{R}{2L}t}\sum_{n=1}^{\infty}\Biggl\{\left(\mathcal{\alpha}_{n}^{2}sin\left(\omega_{n}^{2}t\right)+\beta_{n}^{2}cos\left(\omega_{n}^{2}t\right)\right)\\
\left(\zeta_{n}^{2}sin\left(\sqrt{CL\left(\omega_{n}^{2}\right)^{2}-\frac{CR^{2}}{4L}}x\right)+\eta_{n}^{2}cos\left(\sqrt{CL\left(\omega_{n}^{2}\right)^{2}-\frac{CR^{2}}{4L}}x\right)\right)\Biggr\}.
\end{array}\label{eq: Solution2}
\end{equation}
Here, the superscripts indicated with ``1'' belong to the form of
the solution before the closure of the aortic valve, and the superscripts
indicated with ``2'' belong to the form of the solution after the
closure of the aortic valve. Constants $K^{1}$, $\alpha{}_{n}^{1}$,
$\beta{}_{n}^{1}$, $\zeta{}_{n}^{1}$, $\eta{}_{n}^{1}$ and $\omega{}_{n}^{1}$
are found from the boundary and initial conditions at systole. Similarly,
constants $K^{2}$, $\alpha{}_{n}^{2}$, $\beta{}_{n}^{2}$, $\zeta{}_{n}^{2}$,
$\eta{}_{n}^{2}$ and $\omega{}_{n}^{2}$ are found from the boundary
and initial conditions at diastole. 

Equation (\ref{eq: Solution2}) is explicitly showing the coupling
and decoupling of heart and aorta before and after the dicrotic notch.
As the boundary conditions change during a cardiac cycle, the frequencies
of oscillation also change from $\omega_{n}^{1}$ to $\omega_{n}^{2}$.\textcolor{black}{{}
Generally, Equation (\ref{eq: Solution1}) can represent pressure
waveform for a Newtonian incompressible fluid in a straight and sufficiently
long elastic tube with constant circular cross section. }

If the pressure is recorded at a specific point $x_{0}$ on aorta,
the terms containing the spacial variable $x$ would be fixed. In
other words, Equation (\ref{eq: Solution2}) would reduce to
\begin{equation}
\begin{array}{c}
P\left(x=x_{0},t\right)=\bar{p}\\
+\left\{ K^{1}e^{-\frac{R}{2L}t}\sum_{n=1}^{\infty}\kappa_{n}^{1}\left(\mathcal{\alpha}_{n}^{1}sin\left(\omega_{n}^{1}t\right)+\beta_{n}^{1}cos\left(\omega_{n}^{1}t\right)\right)\right\} \mathbf{1}_{\left[0,T_{0}\right)}\left(t\right)\,\\
+\left\{ K^{2}e^{-\frac{R}{2L}t}\sum_{n=1}^{\infty}\kappa_{n}^{2}\left(\mathcal{\alpha}_{n}^{2}sin\left(\omega_{n}^{2}t\right)+\beta_{n}^{2}cos\left(\omega_{n}^{2}t\right)\right)\right\} \mathbf{1}_{\left[T_{0},T\right)}\left(t\right),
\end{array}\label{eq: Solution at a point}
\end{equation}
for 
\[
\kappa_{n}^{1}=\zeta_{n}^{1}sin\left(\sqrt{CL\left(\omega_{n}^{1}\right)^{2}-\frac{CR^{2}}{4L}}x_{0}\right)+\eta_{n}^{1}cos\left(\sqrt{CL\left(\omega_{n}^{1}\right)^{2}-\frac{CR^{2}}{4L}}x_{0}\right)
\]
 and 
\[
\kappa_{n}^{2}=\zeta_{n}^{2}sin\left(\sqrt{CL\left(\omega_{n}^{2}\right)^{2}-\frac{CR^{2}}{4L}}x_{0}\right)+\eta_{n}^{2}cos\left(\sqrt{CL\left(\omega_{n}^{2}\right)^{2}-\frac{CR^{2}}{4L}}x_{0}\right).
\]
Now, considering that the cardiac cycle length would be around $1.5\,sec$,
at most, and taking into account that $R$ is smaller than $L$ \cite{berger1994differential},
one can use the approximation $e^{-\frac{R}{2L}t}\simeq1$. Hence,
Equation (\ref{eq: Solution at a point}) would become 
\begin{equation}
\begin{array}{c}
P\left(x=x_{0},t\right)\approx\bar{p}\\
+\left\{ K^{1}\sum_{n=1}^{\infty}\kappa_{n}^{1}\left(\mathcal{\alpha}_{n}^{1}sin\left(\omega_{n}^{1}t\right)+\beta_{n}^{1}cos\left(\omega_{n}^{1}t\right)\right)\right\} \mathbf{1}_{\left[0,T_{0}\right)}\left(t\right)\,\\
+\left\{ K^{2}\sum_{n=1}^{\infty}\kappa_{n}^{2}\left(\mathcal{\alpha}_{n}^{2}sin\left(\omega_{n}^{2}t\right)+\beta_{n}^{2}cos\left(\omega_{n}^{2}t\right)\right)\right\} \mathbf{1}_{\left[T_{0},T\right)}\left(t\right).
\end{array}\label{eq: Sol1PointApp}
\end{equation}
Further, if most of the information, or energy, is carried out by
the first terms in the series of the solution, we can further write
the approximated solution (\ref{eq: Sol1PointApp}) as
\begin{equation}
\begin{array}{c}
P\left(x=x_{0},t\right)\approx\bar{p}\\
+\left\{ K^{1}\kappa_{1}^{1}\left(\mathcal{\alpha}_{1}^{1}sin\left(\omega_{1}^{1}t\right)+\beta_{1}^{1}cos\left(\omega_{1}^{1}t\right)\right)\right\} \mathbf{1}_{\left[0,T_{0}\right)}\left(t\right)\,\\
+\left\{ K^{2}\kappa_{1}^{2}\left(\mathcal{\alpha}_{1}^{2}sin\left(\omega_{1}^{2}t\right)+\beta_{1}^{2}cos\left(\omega_{1}^{2}t\right)\right)\right\} \mathbf{1}_{\left[T_{0},T\right)}\left(t\right).
\end{array}\label{eq: OneStepB4IF}
\end{equation}
Now, by relabeling 
\begin{equation}
b_{1}=K^{1}\kappa_{1}^{1}\mathcal{\alpha}_{1}^{1},
\end{equation}
\begin{equation}
a_{1}=K^{1}\kappa_{1}^{1}\beta_{1}^{1},
\end{equation}
\begin{equation}
b_{2}=K^{2}\kappa_{1}^{2}\mathcal{\alpha}_{1}^{2},
\end{equation}
\begin{equation}
a_{2}=K^{2}\kappa_{1}^{2}\beta_{1}^{2},
\end{equation}
\begin{equation}
\omega_{1}=\omega_{1}^{1},
\end{equation}
\begin{equation}
\omega_{2}=\omega_{1}^{2},
\end{equation}
we can approximate the IF model (\ref{eq: The Trend}). The continuity
and periodicity conditions (\ref{eq: Conditions}) can also be approximated
if we hold the assumption that most of the energy is carried out by
the first terms in the series of the solution (\ref{eq: Sol1PointApp}).

In short, in this section, we have presented an approximate quantitative
justification on the origins of the IF method. In the next section,
we move on with the analysis of the optimization problem (\ref{eq: Continuous min})
subject to (\ref{eq: Conditions}). 

\section{Analysis of The IF Algorithm}

Practically, one must solve the discrete version of (\ref{eq: Continuous min}).
We assume that the pressure waveform $f\left(t\right)$ is sampled
uniformly. Also, we can simplify (\ref{eq: Continuous min}) by the
fact that any sinusoid can be assumed to start from time $t=0$ with
a compensation coming from a phase shift. In other words, any sinusoid
can be expressed as $A\cos\omega t+B\sin\omega t$, irrespective of
whether the initial time is $t=0$ or $t=T_{0}$. Hence, the discrete
format of (\ref{eq: Continuous min}) can be expressed as 
\begin{equation}
\begin{array}{cc}
\underset{a_{i},b_{i},\omega_{i},\bar{p}}{minimize} & \left\Vert \mathbf{f}-\mathbf{S}\left(a_{i},b_{i},\bar{p},\omega_{i};\mathbf{t}\right)\right\Vert _{2}^{2}\\
\\
subject\,to & \begin{array}{ccc}
a_{1}\cos\omega_{1}T_{0}+b_{1}\sin\omega_{1}T_{0} & = & a_{2}\phantom{\cos\omega_{2}\left(T-T_{0}\right)+b_{2}\sin\omega_{2}\left(T-T_{0}\right)},\\
a_{1}\phantom{\cos\omega_{1}T_{0}+b_{1}\sin\omega_{1}T_{0}} & = & a_{2}\cos\omega_{2}\left(T-T_{0}\right)+b_{2}\sin\omega_{2}\left(T-T_{0}\right),
\end{array}
\end{array}\label{eq: Discrete min}
\end{equation}
for $\mathbf{f}=\left(f_{1},\ldots,f_{n+m}\right)^{\prime}$ as the
uniform sampling of the original cycle. Here, by taking 
\begin{equation}
\mathbf{t}=\left(\mathbf{t}_{1}^{\prime},\mathbf{t}_{2}^{\prime}\right)^{\prime}=\left(t_{1}^{1},t_{1}^{2},\ldots,t_{1}^{n},t_{2}^{1},t_{2}^{2},\ldots,t_{2}^{m}\right)^{\prime}\in\mathbb{R}^{\left(n+m\right)\times1}\label{eq: times from zero vector}
\end{equation}
for $\mathbf{t}_{1}=\left(0,\Delta t,2\Delta t,\ldots,T_{0}\right)^{\prime}\in\mathbb{R}^{n\times1}$
and $\mathbf{t}_{2}=\left(\Delta t,2\Delta t,\ldots,T-T_{0}\right)^{\prime}\in\mathbb{R}^{m\times1}$,
we have the discrete form of $S\left(a_{i},b_{i},\bar{p},\omega_{i};t\right)$
as
\begin{equation}
\mathbf{S}\left(a_{i},b_{i},\bar{p},\omega_{i};\mathbf{t}\right)=\left(\begin{array}{c}
a_{1}\cos\omega_{1}\mathbf{t}_{1}+b_{1}\sin\omega_{1}\mathbf{t}_{1}\\
a_{2}\cos\omega_{2}\mathbf{t}_{2}+b_{2}\sin\omega_{2}\mathbf{t}_{2}
\end{array}\right)+\bar{p}\mathbf{1}.\label{eq: S form}
\end{equation}
In this article, $\left(.\right)^{\prime}$ denotes the transpose
operator and the vector $\mathbf{1}=\left(1,1,\ldots,1\right)^{\prime}\in\mathbb{R}^{\left(n+m\right)\times1}$.
Also, 
\begin{equation}
\begin{array}{c}
\cos\omega_{1}\mathbf{t}_{1}=\left(\cos\omega_{1}t_{1}^{1},\ldots,\cos\omega_{1}t_{1}^{n}\right)^{\prime},\\
\sin\omega_{1}\mathbf{t}_{1}=\left(\sin\omega_{1}t_{1}^{1},\ldots,\sin\omega_{1}t_{1}^{n}\right)^{\prime},\\
\cos\omega_{2}\mathbf{t}_{2}=\left(\cos\omega_{2}t_{2}^{1},\ldots,\cos\omega_{2}t_{2}^{m}\right)^{\prime},\\
\sin\omega_{2}\mathbf{t}_{2}=\left(\sin\omega_{2}t_{2}^{1},\ldots,\sin\omega_{2}t_{2}^{m}\right)^{\prime}.
\end{array}
\end{equation}
The constraints, in (\ref{eq: Discrete min}), can be written as 
\begin{equation}
\left(\begin{array}{cccc}
\cos\omega_{1}T_{0} & -1 & \sin\omega_{1}T_{0} & 0\\
1 & -\cos\omega_{2}\left(T-T_{0}\right) & 0 & -\sin\omega_{2}\left(T-T_{0}\right)
\end{array}\right)\left(\begin{array}{c}
a_{1}\\
a_{2}\\
b_{1}\\
b_{2}
\end{array}\right)=\left(\begin{array}{c}
0\\
0
\end{array}\right).\label{eq: constraint linear eq}
\end{equation}
If we can solve for two, out of four, unknowns in (\ref{eq: constraint linear eq}),
we would make (\ref{eq: Discrete min}) an unconstrained optimization.
However, it is important to check whether the matrix in (\ref{eq: constraint linear eq})
is of full rank or not. In fact, the rows of this matrix are linearly
independent except when 
\begin{equation}
\cos\omega_{1}T_{0}\cos\omega_{2}\left(T-T_{0}\right)=1.\label{eq: Node Locations}
\end{equation}
This will lead into two cases: 
\begin{enumerate}
\item \emph{Degenerate Case} in which Equation (\ref{eq: Node Locations})
holds,
\item \emph{General Case} in which, it does not. 
\end{enumerate}

\subsection{General Case ($\cos\omega_{1}T_{0}\cos\omega_{2}\left(T-T_{0}\right)\protect\neq1$)}

One can solve the constraints in (\ref{eq: Discrete min}) for $a_{1}$
and $a_{2}$ to obtain
\begin{equation}
a_{1}=\frac{b_{1}\sin\omega_{1}T_{0}\cos\omega_{2}\left(T-T_{0}\right)+b_{2}\sin\omega_{2}\left(T-T_{0}\right)}{1-\cos\omega_{1}T_{0}\cos\omega_{2}\left(T-T_{0}\right)},\label{eq: a1}
\end{equation}
\begin{equation}
a_{2}=\frac{b_{1}\sin\omega_{1}T_{0}+b_{2}\cos\omega_{1}T_{0}\sin\omega_{2}\left(T-T_{0}\right)}{1-\cos\omega_{1}T_{0}\cos\omega_{2}\left(T-T_{0}\right)}.\label{eq: a2}
\end{equation}
Equations (\ref{eq: a1}) and (\ref{eq: a2}) would then simplify
(\ref{eq: S form}) into 
\begin{equation}
\mathbf{S}\left(\omega_{1},\omega_{1},b_{1},b_{2},\bar{p};\mathbf{t}\right)=\mathbf{Q}\left(\omega_{1},\omega_{1},b_{1},b_{2};\mathbf{t}\right)+\bar{p}\mathbf{1},\label{eq: Simplified S}
\end{equation}
where $\mathbf{Q}\left(\omega_{1},\omega_{1},b_{1},b_{2};\mathbf{t}\right)=b_{1}\mathbf{v}_{1}\left(\omega_{1},\omega_{2};\mathbf{t}\right)+b_{2}\mathbf{v}_{2}\left(\omega_{1},\omega_{2};\mathbf{t}\right)$
for 
\begin{equation}
\mathbf{v}_{1}\left(\omega_{1},\omega_{2};\mathbf{t}\right)=\left(\begin{array}{c}
\frac{\sin\omega_{1}T_{0}\cos\omega_{2}\left(T-T_{0}\right)}{1-\cos\omega_{1}T_{0}\cos\omega_{2}\left(T-T_{0}\right)}\cos\omega_{1}\mathbf{t}_{1}+\sin\omega_{1}\mathbf{t}_{1}\\
\frac{\sin\omega_{1}T_{0}}{1-\cos\omega_{1}T_{0}\cos\omega_{2}\left(T-T_{0}\right)}\cos\omega_{2}\mathbf{t}_{2}\phantom{+\sin\omega_{1}\mathbf{t}_{1}}
\end{array}\right),\label{eq: v1}
\end{equation}
and
\begin{equation}
\mathbf{v}_{2}\left(\omega_{1},\omega_{2};\mathbf{t}\right)=\left(\begin{array}{c}
\frac{\sin\omega_{2}\left(T-T_{0}\right)}{1-\cos\omega_{1}T_{0}\cos\omega_{2}\left(T-T_{0}\right)}\cos\omega_{1}\mathbf{t}_{1}\phantom{+\sin\omega_{2}\mathbf{t}_{2}}\\
\frac{\cos\omega_{1}T_{0}\sin\omega_{2}\left(T-T_{0}\right)}{1-\cos\omega_{1}T_{0}\cos\omega_{2}\left(T-T_{0}\right)}\cos\omega_{2}\mathbf{t}_{2}+\sin\omega_{2}\mathbf{t}_{2}
\end{array}\right).\label{eq: v2}
\end{equation}
Using Equations (\ref{eq: Simplified S})-(\ref{eq: v2}), and dropping
the dependencies in notation, simplifies (\ref{eq: Discrete min})
into 
\begin{equation}
\underset{\omega_{1},\omega_{2},b_{1},b_{2},\bar{p}}{minimize}\,\left\Vert \mathbf{Q}+\bar{p}\mathbf{1}-\mathbf{f}\right\Vert _{2}^{2}.\label{eq: No-constraint Minimization}
\end{equation}
This simplification has helped to eliminate the constraints in the
optimization problem (\ref{eq: Discrete min}). 

The minimization problem (\ref{eq: No-constraint Minimization}) is
non-convex and non-linear in its parameters. So, in order to be able
to solve the problem, we can use the fact that the minimum of a function
can first be found over some variables and then over the remaining
ones \cite{boyd2004convex}. In other words, the optimization problem
in (\ref{eq: No-constraint Minimization}) can be written as 
\begin{equation}
\underset{\omega_{1},\omega_{2}}{minimize}\,\left(\underset{b_{1},b_{2},\bar{p}}{minimize}\,\left\Vert \mathbf{Q}+\bar{p}\mathbf{1}-\mathbf{f}\right\Vert _{2}^{2}\right).\label{eq: Two Step Minimization}
\end{equation}

We call the inner optimization in (\ref{eq: Two Step Minimization})
as $P\left(\omega_{1},\omega_{2}\right)$. Solving for $P\left(\omega_{1},\omega_{2}\right)$
is a classical least squares problem. The solution existence and uniqueness
of this optimization is mentioned in our previous work \cite{tavallali2015convergence}.
To find the exact solution we simplify the objective function as 

\begin{equation}
\begin{array}{cc}
\left\Vert \mathbf{Q}+\bar{p}\mathbf{1}-\mathbf{f}\right\Vert _{2}^{2}= & \left(\mathbf{Q}+\bar{p}\mathbf{1}-\mathbf{f}\right)^{\prime}\left(\mathbf{Q}+\bar{p}\mathbf{1}-\mathbf{f}\right)\\
\phantom{\left\Vert \mathbf{Q}+\bar{p}\mathbf{1}-\mathbf{f}\right\Vert _{2}^{2}}= & \mathbf{Q}^{\prime}\mathbf{Q}+2\bar{p}\mathbf{Q}^{\prime}\mathbf{1}-2\mathbf{Q}^{\prime}\mathbf{f}-2\bar{p}\mathbf{f}^{\prime}\mathbf{1}+\bar{p}\mathbf{1}^{\prime}\mathbf{1}+\mathbf{f}^{\prime}\mathbf{f}.
\end{array}\label{eq: Norm to dot product}
\end{equation}
Substituting for $\mathbf{Q}=b_{1}\mathbf{v}_{1}+b_{2}\mathbf{v}_{2}$,
we convert (\ref{eq: Norm to dot product}) into 
\begin{equation}
\begin{array}{cc}
\left\Vert \mathbf{Q}+\bar{p}\mathbf{1}-\mathbf{f}\right\Vert _{2}^{2}= & b_{1}^{2}\mathbf{v}_{1}^{\prime}\mathbf{v}_{1}+2b_{1}b_{2}\mathbf{v}_{1}^{\prime}\mathbf{v}_{2}+b_{2}^{2}\mathbf{v}_{2}^{\prime}\mathbf{v}_{2}+2\bar{p}b_{1}\mathbf{v}_{1}^{\prime}\mathbf{1}+2\bar{p}b_{2}\mathbf{v}_{2}^{\prime}\mathbf{1}\\
\phantom{\left\Vert \mathbf{Q}+\bar{p}\mathbf{1}-\mathbf{f}\right\Vert _{2}^{2}=} & -2b_{1}\mathbf{v}_{1}^{\prime}\mathbf{f}-2b_{2}\mathbf{v}_{2}^{\prime}\mathbf{f}-2\bar{p}\mathbf{f}^{\prime}\mathbf{1}+\bar{p}^{2}\mathbf{1}^{\prime}\mathbf{1}+\mathbf{f}^{\prime}\mathbf{f}.
\end{array}\label{eq: Ready for Derivative}
\end{equation}
Since, in this part of the optimization, the values of $\omega_{1}$
and $\omega_{2}$ are fixed, we can find the optimal values of $b_{1}$,
$b_{2}$, and $\bar{p}$ by setting the partial derivatives of (\ref{eq: Ready for Derivative})
equal to zero. In other words, we set $\frac{\partial\left(\left\Vert \mathbf{Q}+\bar{p}\mathbf{1}-\mathbf{f}\right\Vert _{2}^{2}\right)}{\partial b_{1}}=0$,
$\frac{\partial\left(\left\Vert \mathbf{Q}+\bar{p}\mathbf{1}-\mathbf{f}\right\Vert _{2}^{2}\right)}{\partial b_{2}}=0$,
and $\frac{\partial\left(\left\Vert \mathbf{Q}+\bar{p}\mathbf{1}-\mathbf{f}\right\Vert _{2}^{2}\right)}{\partial\bar{p}}=0$.
Doing this, we find the optimal solution for $b_{1}$, $b_{2}$, and
$\bar{p}$, by 
\begin{equation}
\left(\begin{array}{c}
b_{1}^{*}\left(\omega_{1},\omega_{2}\right)\\
b_{2}^{*}\left(\omega_{1},\omega_{2}\right)\\
\bar{p}^{*}\left(\omega_{1},\omega_{2}\right)
\end{array}\right)=\left(\begin{array}{ccc}
\mathbf{v}_{1}^{\prime}\mathbf{v}_{1} & \mathbf{v}_{1}^{\prime}\mathbf{v}_{2} & \mathbf{v}_{1}^{\prime}\mathbf{1}\\
\mathbf{v}_{1}^{\prime}\mathbf{v}_{2} & \mathbf{v}_{2}^{\prime}\mathbf{v}_{2} & \mathbf{v}_{2}^{\prime}\mathbf{1}\\
\mathbf{v}_{1}^{\prime}\mathbf{1} & \mathbf{v}_{2}^{\prime}\mathbf{1} & \mathbf{1}^{\prime}\mathbf{1}
\end{array}\right)^{-1}\left(\begin{array}{c}
\mathbf{v}_{1}^{\prime}\mathbf{f}\\
\mathbf{v}_{2}^{\prime}\mathbf{f}\\
\mathbf{1}^{\prime}\mathbf{f}
\end{array}\right).\label{eq: Inner loop argmin sol}
\end{equation}
Here, we have fulfilled the optimization part by solving a linear
system. This could potentially accelerate the IF algorithm. Finally,
we only have to solve a minimization on 
\begin{equation}
P\left(\omega_{1},\omega_{2}\right)=\left\Vert \mathbf{Q}\left(\omega_{1},\omega_{2},b_{1}^{*}\left(\omega_{1},\omega_{2}\right),b_{2}^{*}\left(\omega_{1},\omega_{2}\right);\mathbf{t}\right)+\bar{p}^{*}\left(\omega_{1},\omega_{2}\right)\mathbf{1}-\mathbf{f}\right\Vert _{2}^{2},\label{eq: P definition}
\end{equation}
which is 
\begin{equation}
\underset{\omega_{1},\omega_{2}}{minimize}\,P\left(\omega_{1},\omega_{2}\right).\label{eq: P minimization}
\end{equation}
We note that a property of the function $P\left(\omega_{1},\omega_{2}\right)$
is its differentiability, away from its singularities. In fact, by
definition, the function $\left\Vert \mathbf{Q}+\bar{p}\mathbf{1}-\mathbf{f}\right\Vert _{2}^{2}$
is directionally differentiable with respect to all its variables.
Hence, using the results in \cite{borisenko1992directional,pshenichnyi1971necessary},
we can deduce that 
\begin{equation}
P\left(\omega_{1},\omega_{2}\right)=\underset{b_{1},b_{2},\bar{p}}{minimize}\,\left\Vert \mathbf{Q}+\bar{p}\mathbf{1}-\mathbf{f}\right\Vert _{2}^{2}
\end{equation}

is directionally differentiable with respect to $\omega_{1}$ and
$\omega_{2}$. This property can be exploited if one tries to solve
(\ref{eq: P minimization}) using a gradient based optimization method
\cite{bertsekas1999nonlinear}. 

\subsection{Degenerate Case ($\cos\omega_{1}T_{0}\cos\omega_{2}\left(T-T_{0}\right)=1$)}

The solution of (\ref{eq: Node Locations}) can be expressed as nodes
of a lattice $\mathcal{N}$ in $\omega_{1}\omega_{2}$ plane. To be
more specific, we have
\begin{equation}
\mathcal{N}=\varGamma_{1}\cup\varGamma_{2},\label{eq: Nodes}
\end{equation}
where
\begin{equation}
\varGamma_{1}=\left\{ \left(\omega_{1},\omega_{2}\right)\left|\omega_{1}T_{0}=\left(2k_{1}+1\right)\pi,\omega_{2}\left(T-T_{0}\right)=\left(2k_{2}+1\right)\pi,k_{1}\in\mathbb{Z},k_{2}\in\mathbb{Z}\right.\right\} ,\label{eq: Lattice Nodes1}
\end{equation}
and 
\begin{equation}
\begin{array}{c}
\varGamma_{2}=\left\{ \left(\omega_{1},\omega_{2}\right)\left|\omega_{1}T_{0}=2k_{1}\pi,\omega_{2}\left(T-T_{0}\right)=2k_{2}\pi,k_{1}\in\mathbb{Z},k_{2}\in\mathbb{Z}\right.\right\} .\end{array}\label{eq: Lattice Nodes2}
\end{equation}
If $\left(\omega_{1},\omega_{2}\right)\in\varGamma_{1}$, from (\ref{eq: constraint linear eq})
we have $a_{1}=-a_{2}$. On the other hand, if $\left(\omega_{1},\omega_{2}\right)\in\varGamma_{2}$,
from (\ref{eq: constraint linear eq}) we have $a_{1}=a_{2}$. In
both of these cases, we can express (\ref{eq: S form}) as 
\begin{equation}
\mathbf{S}\left(\omega_{1},\omega_{1},a_{1},b_{1},b_{2},\bar{p};\mathbf{t}\right)=\mathbf{Q}\left(\omega_{1},\omega_{1},a_{1},b_{1},b_{2};\mathbf{t}\right)+\bar{p}\mathbf{1},\label{eq: Simplified S Degen}
\end{equation}
where $\mathbf{Q}\left(\omega_{1},\omega_{1},a_{1},b_{1},b_{2};\mathbf{t}\right)=a_{1}\mathbf{w}_{0}^{\varGamma_{i}}\left(\omega_{1},\omega_{2};\mathbf{t}\right)+b_{1}\mathbf{w}_{1}\left(\omega_{1},\omega_{2};\mathbf{t}\right)+b_{2}\mathbf{w}_{2}\left(\omega_{1},\omega_{2};\mathbf{t}\right)$,
for $i=1,2$. If $\left(\omega_{1},\omega_{2}\right)\in\varGamma_{1}$,
\begin{equation}
\mathbf{w}_{0}^{\varGamma_{1}}=\left(\begin{array}{c}
\begin{array}{c}
\phantom{-}\cos\omega_{1}\mathbf{t}_{1}\end{array}\\
-\cos\omega_{2}\mathbf{t}_{2}
\end{array}\right).
\end{equation}
Similarly, if $\left(\omega_{1},\omega_{2}\right)\in\varGamma_{2}$,
we have 
\begin{equation}
\mathbf{w}_{0}^{\varGamma_{2}}=\left(\begin{array}{c}
\begin{array}{c}
\cos\omega_{1}\mathbf{t}_{1}\end{array}\\
\cos\omega_{2}\mathbf{t}_{2}
\end{array}\right).
\end{equation}
In both of the cases, we have 
\begin{equation}
\mathbf{w}_{1}=\left(\begin{array}{c}
\sin\omega_{1}\mathbf{t}_{1}\\
\mathbf{0}_{1}
\end{array}\right),\label{eq: w1}
\end{equation}
and 
\begin{equation}
\mathbf{w}_{2}=\left(\begin{array}{c}
\mathbf{0}_{2}\\
\sin\omega_{2}\mathbf{t}_{2}
\end{array}\right).\label{eq: w2}
\end{equation}
Here, $\mathbf{0}_{1}$ and $\mathbf{0}_{2}$ are zero vectors in
$\mathbb{R}^{m\times1}$ and $\mathbb{R}^{n\times1}$, respectively.
It is clear, from (\ref{eq: w1}) and (\ref{eq: w2}), that $\mathbf{w}_{1}^{\prime}\mathbf{w}_{2}=\mathbf{w}_{2}^{\prime}\mathbf{w}_{1}=0$.
Using (\ref{eq: Simplified S Degen}), and a similar approach we employed
in (\ref{eq: Ready for Derivative}) and (\ref{eq: Inner loop argmin sol}),
we find the optimal solution for $a_{1}$, $b_{1}$, $b_{2}$, and
$\bar{p}$, by 
\begin{equation}
\left(\begin{array}{c}
a_{1,i}^{*}\left(\omega_{1},\omega_{2}\right)\\
b_{1,i}^{*}\left(\omega_{1},\omega_{2}\right)\\
b_{2,i}^{*}\left(\omega_{1},\omega_{2}\right)\\
\bar{p}_{,i}^{*}\left(\omega_{1},\omega_{2}\right)
\end{array}\right)=\left(\begin{array}{cccc}
\left(\mathbf{w}_{0}^{\varGamma_{i}}\right)^{\prime}\mathbf{w}_{0}^{\varGamma_{i}} & \left(\mathbf{w}_{0}^{\varGamma_{i}}\right)^{\prime}\mathbf{w}_{1} & \left(\mathbf{w}_{0}^{\varGamma_{i}}\right)^{\prime}\mathbf{w}_{2} & \left(\mathbf{w}_{0}^{\varGamma_{i}}\right)^{\prime}\mathbf{1}\\
\left(\mathbf{w}_{0}^{\varGamma_{i}}\right)^{\prime}\mathbf{w}_{1} & \mathbf{w}_{1}^{\prime}\mathbf{w}_{1} & 0 & \mathbf{w}_{1}^{\prime}\mathbf{1}\\
\left(\mathbf{w}_{0}^{\varGamma_{i}}\right)^{\prime}\mathbf{w}_{2} & 0 & \mathbf{w}_{2}^{\prime}\mathbf{w}_{2} & \mathbf{w}_{2}^{\prime}\mathbf{1}\\
\left(\mathbf{w}_{0}^{\varGamma_{i}}\right)^{\prime}\mathbf{1} & \mathbf{w}_{1}^{\prime}\mathbf{1} & \mathbf{w}_{2}^{\prime}\mathbf{1} & \mathbf{1}^{\prime}\mathbf{1}
\end{array}\right)^{-1}\left(\begin{array}{c}
\left(\mathbf{w}_{0}^{\varGamma_{i}}\right)^{\prime}\mathbf{f}\\
\mathbf{w}_{1}^{\prime}\mathbf{f}\\
\mathbf{w}_{2}^{\prime}\mathbf{f}\\
\mathbf{1}^{\prime}\mathbf{f}
\end{array}\right),
\end{equation}
for $i=1,2$. Hence, similar to (\ref{eq: P definition}), for $\left(\omega_{1},\omega_{2}\right)\in\varGamma_{1}$
or $\left(\omega_{1},\omega_{2}\right)\in\varGamma_{2}$, we only
have to solve a minimization on 
\begin{equation}
P\left(\omega_{1},\omega_{2}\right)=\left\Vert \mathbf{Q}\left(\omega_{1},\omega_{1},a_{1,i}^{*}\left(\omega_{1},\omega_{2}\right),b_{1,i}^{*}\left(\omega_{1},\omega_{2}\right),b_{2,i}^{*}\left(\omega_{1},\omega_{2}\right);\mathbf{t}\right)+\bar{p}_{,i}^{*}\left(\omega_{1},\omega_{2}\right)\mathbf{1}-\mathbf{f}\right\Vert _{2}^{2}.\label{eq: P definition-1}
\end{equation}
Note that, from a machine learning perspective, the nodes specified
in (\ref{eq: Nodes}) do not have important information physiologically
as they could be inferred from the systolic and diastolic parts of
a waveform alone. In other words, even if these points present a global
minima, they are not informative as we already know the systolic and
diastolic inverses, $\frac{1}{T_{0}}$ and $\frac{1}{T-T_{0}}$ respectively,
as possible inputs to any machine learning algorithm. Hence, these
points could possibly be ignored in a search for an optimum point
of (\ref{eq: Discrete min}). 

\section{Fast IF Algorithms}

In this section, we present a fast IF algorithm which is based on
the results presented in the previous section and the topology of
the solution space for $P\left(\omega_{1},\omega_{2}\right)$. In
order to keep the fluency of this section, we mention the original
IF algorithm (see Algorithm \ref{Alg: InF Algorithm}) as presented
in \cite{tavallali2015convergence}. 

Algorithm \ref{Alg: InF Algorithm} has three major steps. In the
first step, the $\left(\omega_{1},\omega_{2}\right)$ domain 
\begin{equation}
\mathcal{D}_{fr}=\left\{ \left(\omega_{1},\omega_{2}\right)\left|0<\omega_{1}\leq C,\,0<\omega_{2}\leq C\right.\right\} \label{eq: Domain}
\end{equation}
is made discrete, namely $\mathcal{\overline{D}}_{fr}$. The second
step is a minimization to find $P\left(\omega_{1},\omega_{2}\right)$,
see (\ref{eq: P minimization}). The final step is a brute-force search
on $\mathcal{\overline{D}}_{fr}$ to find the minimum of $P\left(\omega_{1},\omega_{2}\right)$. 

All three steps can be optimized to make the IF algorithm faster.
Regarding the domain of optimization $\mathcal{D}_{fr}$, defined
in (\ref{eq: Domain}), we know from our previous work in \cite{pahlevan2017noninvasive}
that the average IF solution, for a physiological pulse waveform recording,
is confined to a smaller domain $\mathcal{D}$ expressed as
\begin{equation}
\mathcal{D}=\left\{ \left(\omega_{1},\omega_{2}\right)\left|0.5\leqslant\frac{\omega_{1}T_{0}}{\pi}\leqslant1.5,\,0.5\leqslant\frac{\omega_{2}\left(T-T_{0}\right)}{\pi}\leqslant3\right.\right\} .\label{eq: Better Domain}
\end{equation}
This will make the first step search area more well-defined and optimized.
In the previous section, we have been able to find some analytic solutions
(see (\ref{eq: Inner loop argmin sol})) for the inner optimization
part of problem (\ref{eq: Two Step Minimization}). This will help
us to substitute an analytic solution instead of an iterative \cite{gould2004preprocessing}
or QR decomposition \cite{trefethen1997numerical} solution for (\ref{eq: Algorithm minimization}).
Finally, the brute-force part can be substituted with an appropriate
direct search algorithm \cite{kolda2003optimization}, e.g. pattern
search algorithm \cite{hooke1961direct}. It can even be substituted
with an appropriate gradient based algorithm \cite{bertsekas1999nonlinear,boyd2004convex},
e.g. gradient descent, as we know the differentiability of $P\left(\omega_{1},\omega_{2}\right)$. 

\begin{algorithm}[H]
\begin{enumerate}
\item Make $\mathcal{D}_{fr}$ discrete for a uniform $r\times r$ mesh
$\mathcal{\overline{D}}_{fr}$, $r\in\mathbb{N}$, 
\[
\mathcal{\overline{D}}_{fr}=\left\{ \left(\omega_{1}^{l},\omega_{2}^{m}\right)\left|\omega_{1}=\frac{l}{r}C,\,\omega_{2}=\frac{m}{r}C;\,l,m\in\left\{ 0,1,\ldots,r\right\} \right.\right\} .
\]
\item For all $l,m\in\left\{ 0,1,\ldots,r\right\} $ solve
\begin{equation}
\begin{array}{cc}
\underset{a_{i},b_{i},\bar{p}}{minimize} & \sum_{j=1}^{n}\left(f\left(t_{j}\right)-S\left(a_{i},b_{i},\bar{p},\omega_{1}^{l},\omega_{2}^{m};t_{j}\right)\right)^{2}\\
\\
subject\,to & \begin{array}{ccc}
a_{1}\cos\omega_{1}T_{0}+b_{1}\sin\omega_{1}T_{0} & = & a_{2}\cos\omega_{2}T_{0}+b_{2}\sin\omega_{2}T_{0},\\
a_{1} & = & a_{2}\cos\omega_{2}T+b_{2}\sin\omega_{2}T.
\end{array}
\end{array}\label{eq: Algorithm minimization}
\end{equation}
and store $P\left(\omega_{1}^{l},\omega_{2}^{m}\right)=\sum_{j=1}^{n}\left(f\left(t_{j}\right)-S\left(a_{i}^{*},b_{i}^{*},\bar{p}^{*},\omega_{1}^{l},\omega_{2}^{m};t_{j}\right)\right)^{2}$
for minimizers $a_{i}^{*},b_{i}^{*},\bar{p}^{*}$.
\item Find the intrinsic frequencies (IFs) 
\[
\left(\omega_{1}^{*},\omega_{2}^{*}\right)=\underset{l,m}{argmin}\left(P\left(\omega_{1}^{l},\omega_{2}^{m}\right)\right).
\]
\end{enumerate}
\caption{Intrinsic Frequency}

\label{Alg: InF Algorithm}
\end{algorithm}

\begin{figure}[H]
\includegraphics[scale=0.3]{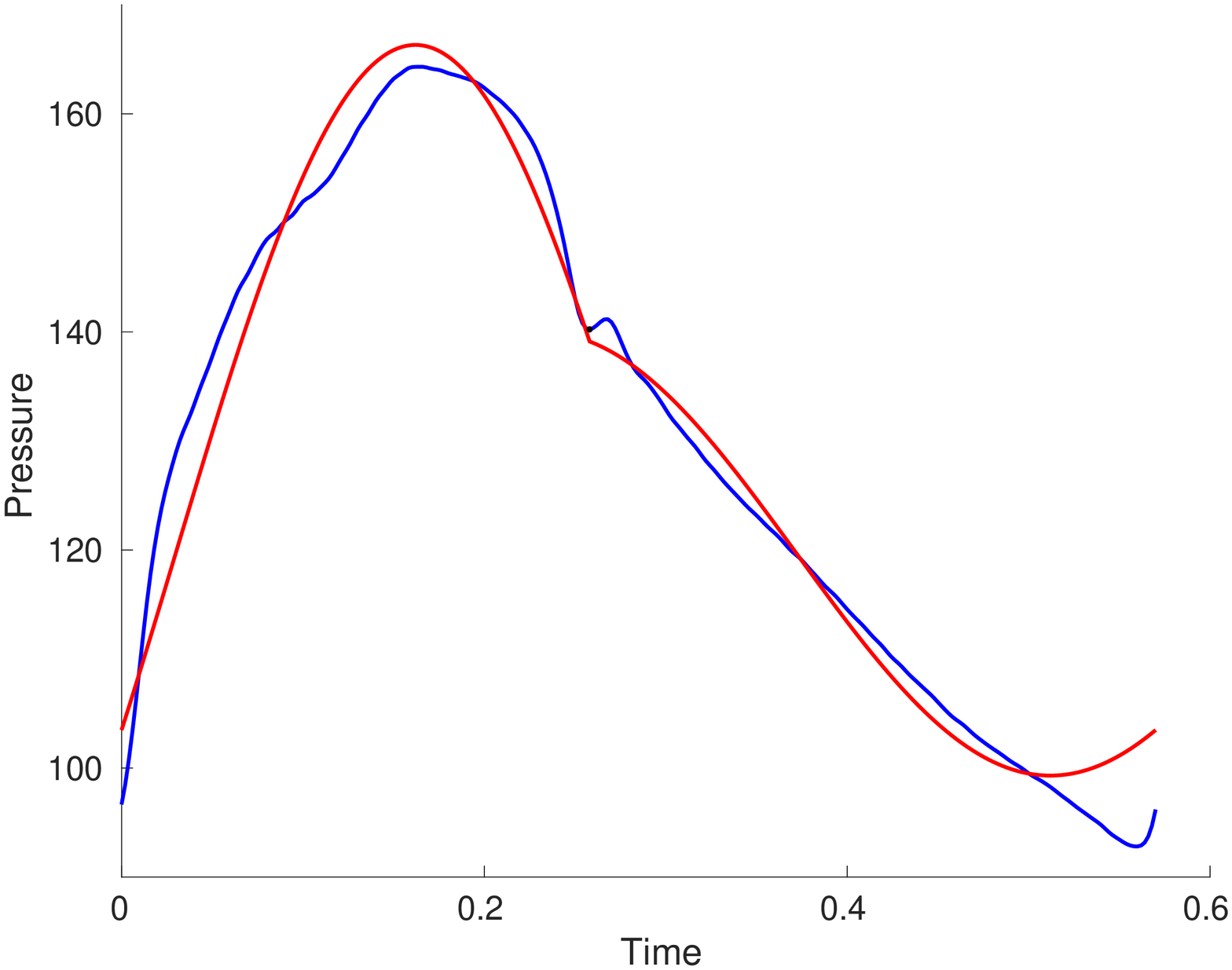}\includegraphics[scale=0.3]{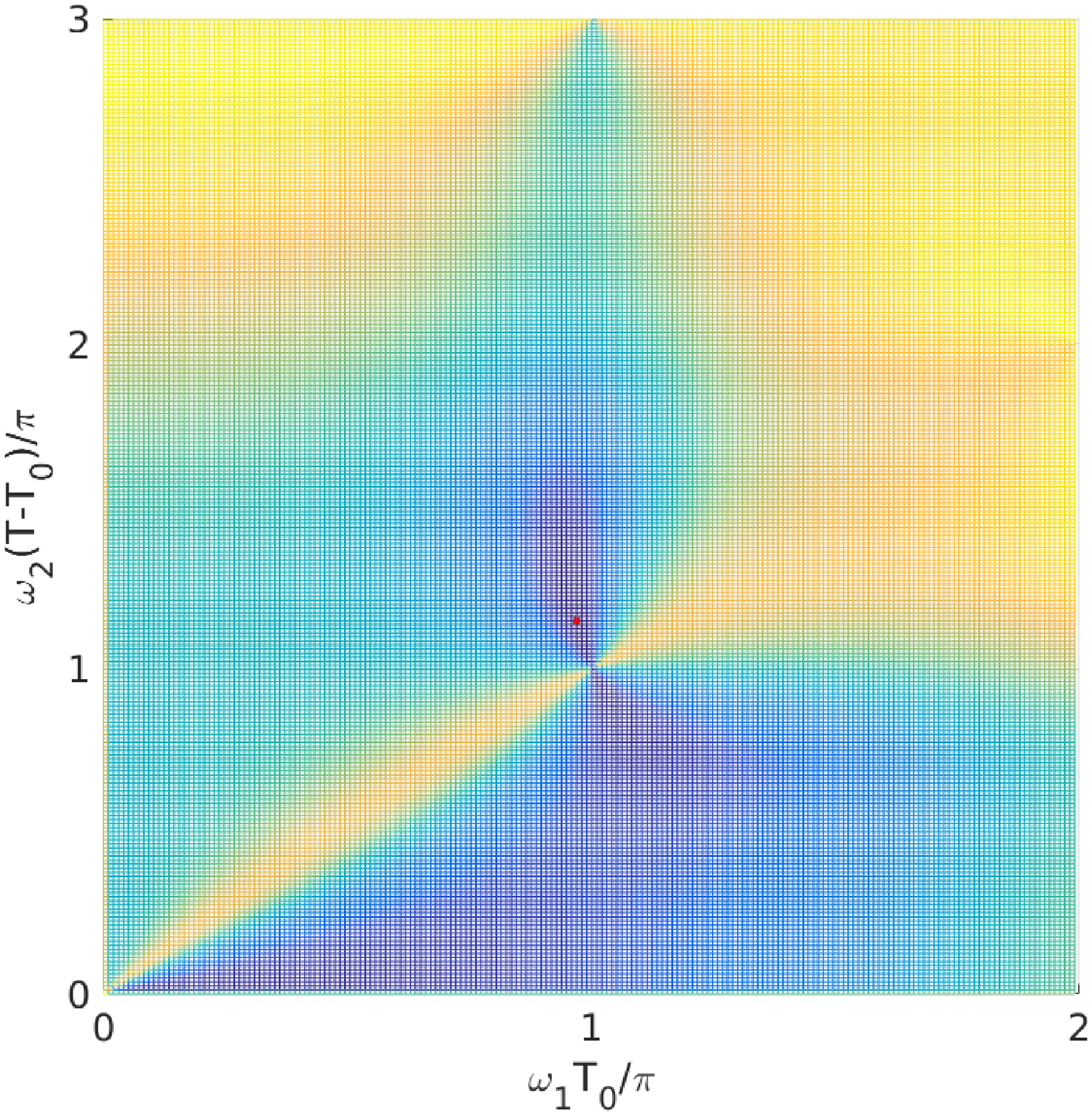}

\includegraphics[scale=0.3]{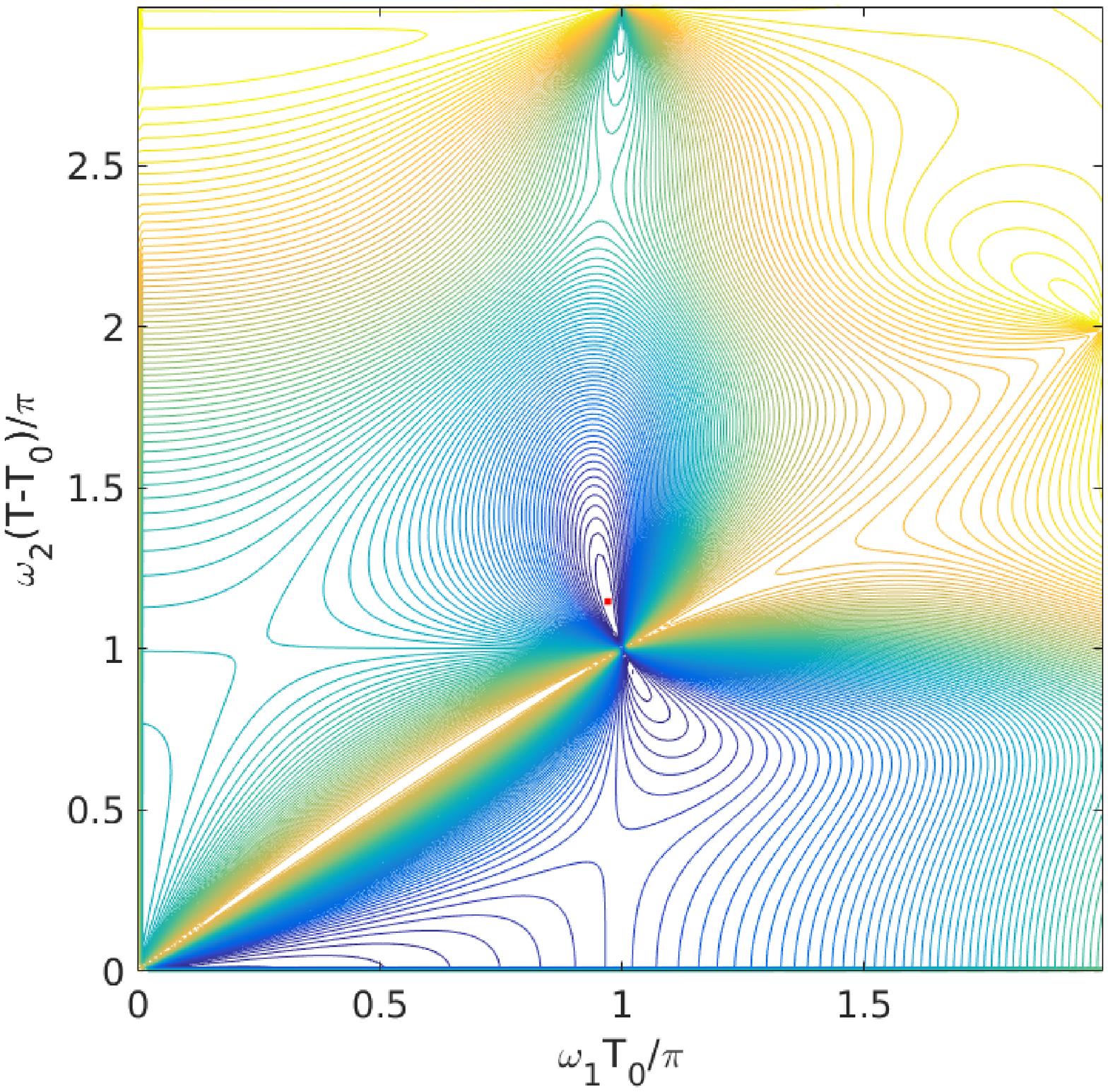}\includegraphics[scale=0.3]{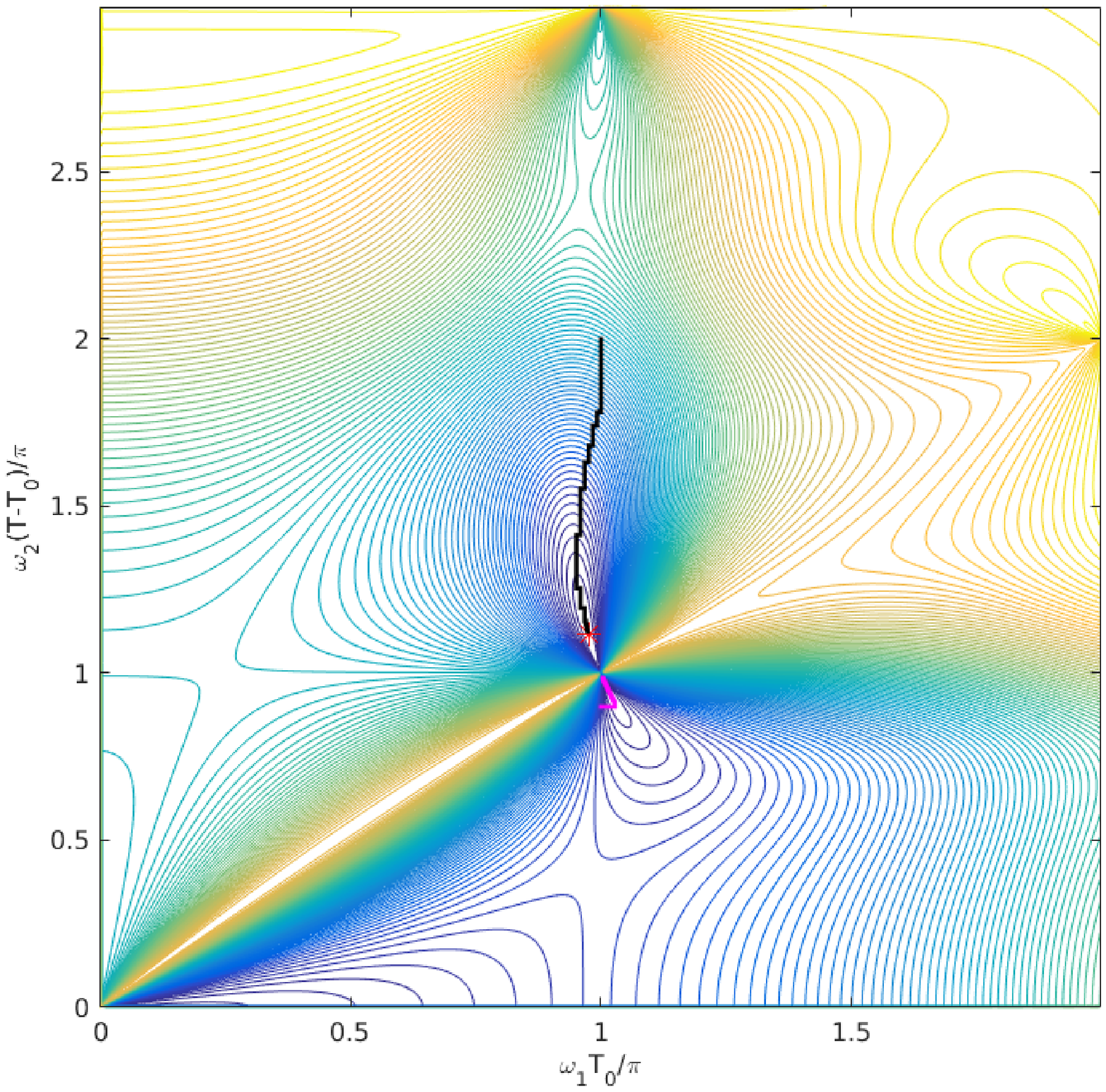}

\caption{\textbf{Up-Left:} A dog aortic pressure cycle (in blue), its dicrotic
notch (black dot), and the IMF (in red). \textbf{Up-Right:} heat-map
plot of $P\left(\frac{\omega_{1}T_{0}}{\pi},\frac{\omega_{2}\left(T-T_{0}\right)}{\pi}\right)$
for the cycle in left with the location of the solution marked with
red dot. \textbf{Down-Left:} Contour plot of $P\left(\frac{\omega_{1}T_{0}}{\pi},\frac{\omega_{2}\left(T-T_{0}\right)}{\pi}\right)$.
The location of the minimizer of $P$ is shown by a red dot. \textbf{Down-Right:}
Contour plot of $P\left(\frac{\omega_{1}T_{0}}{\pi},\frac{\omega_{2}\left(T-T_{0}\right)}{\pi}\right)$
and the location of the minimizer of $P$ tracked by the pattern search.
The true optimum point is marked with a red asterisk. The upper pattern
search set (in black) has converged towards the correct optimum. The
lower pattern search set (in magenta) has converged to a local minima
near the node. }

\label{Fig: Fit plot and Mesh 1}
\end{figure}

\begin{figure}[H]
\includegraphics[scale=0.3]{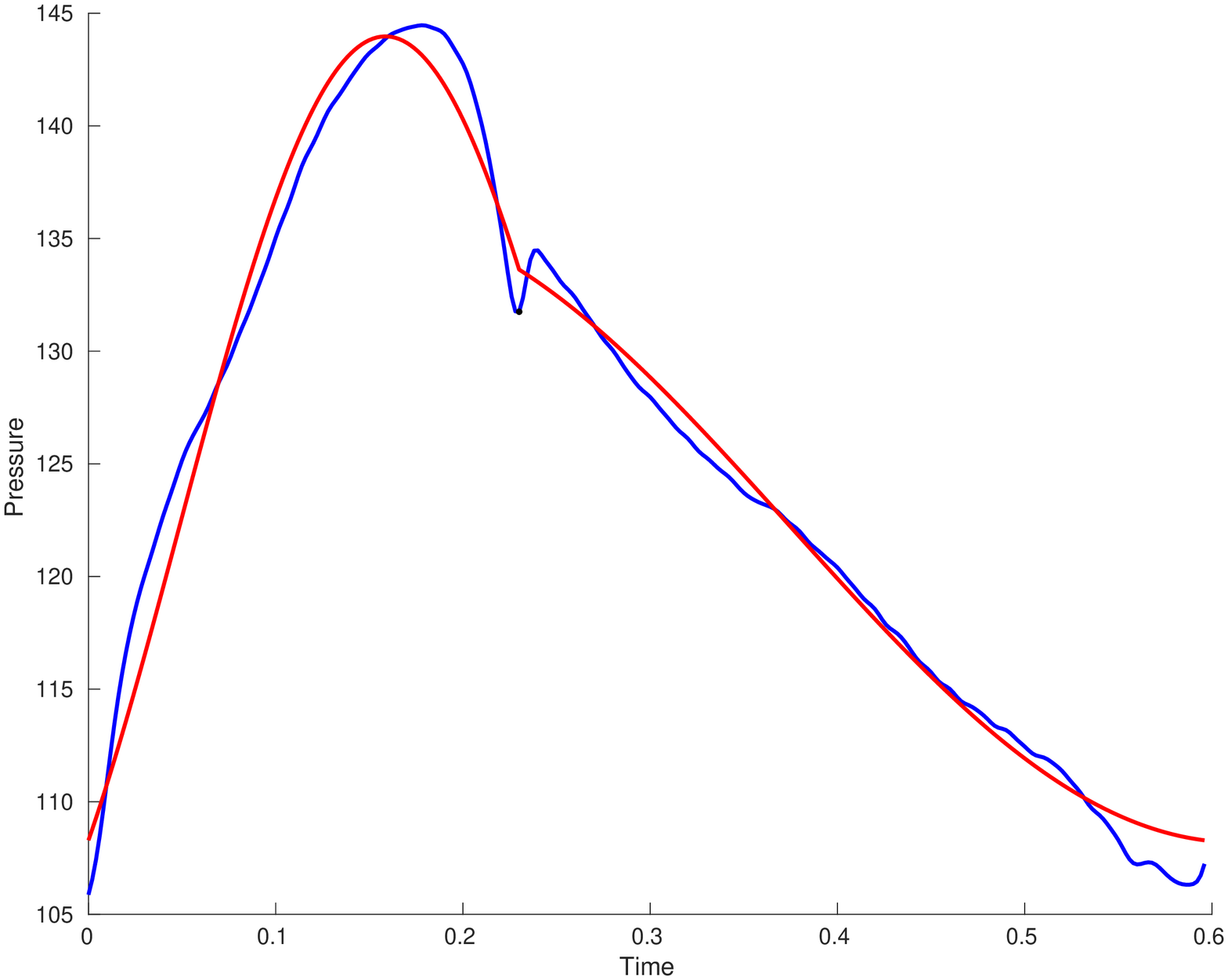}\includegraphics[scale=0.3]{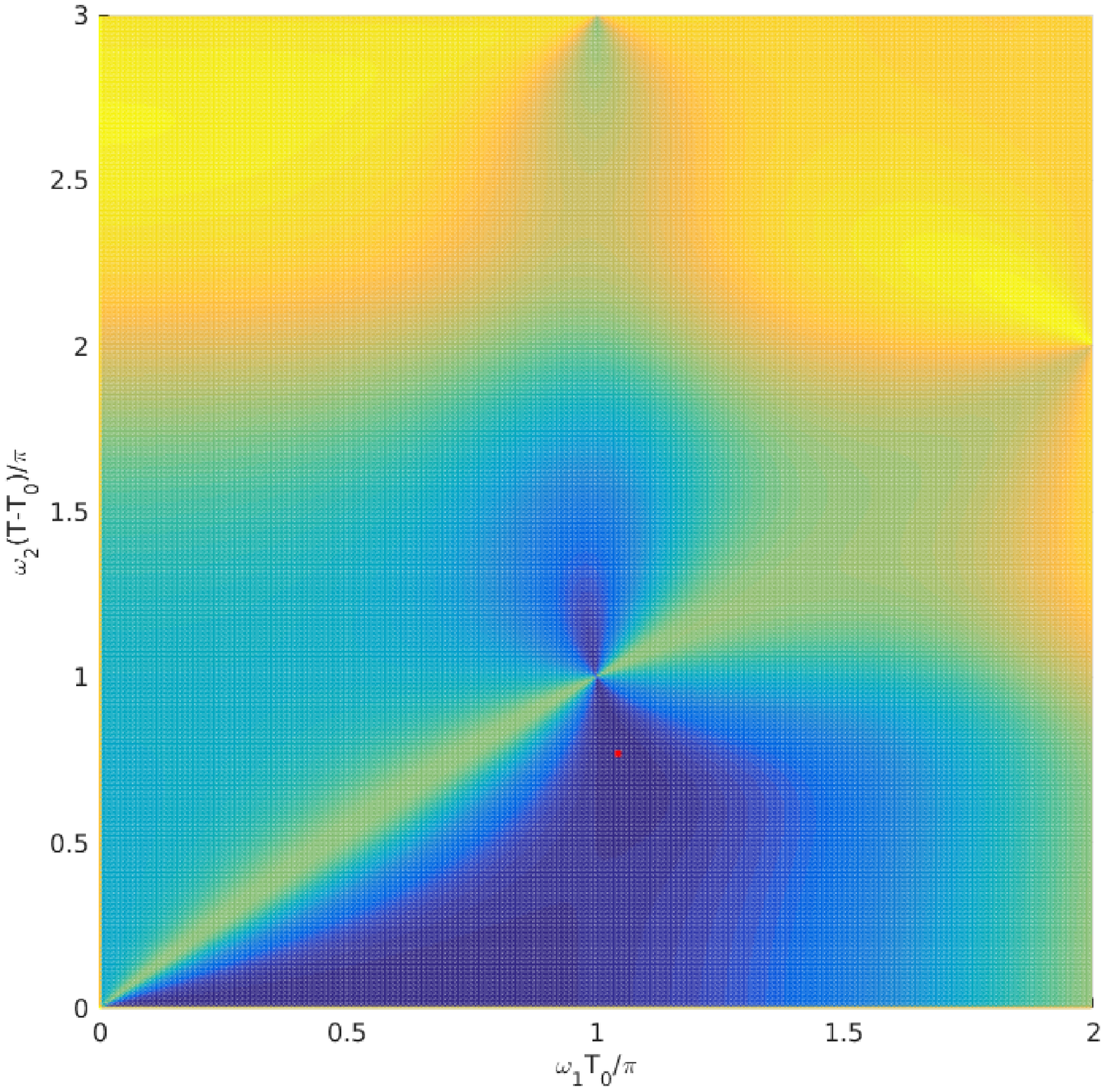}

\includegraphics[scale=0.3]{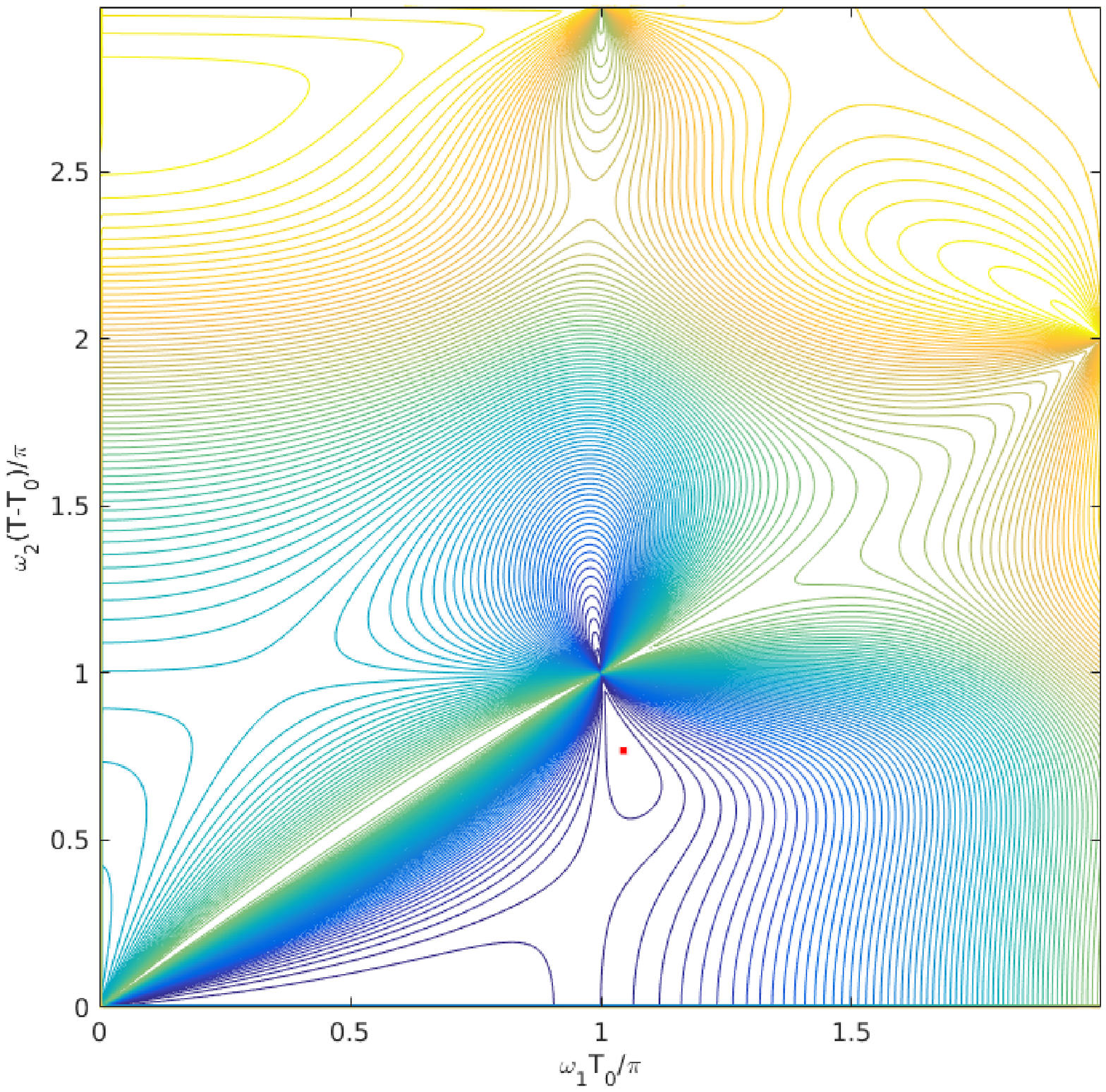}\includegraphics[scale=0.3]{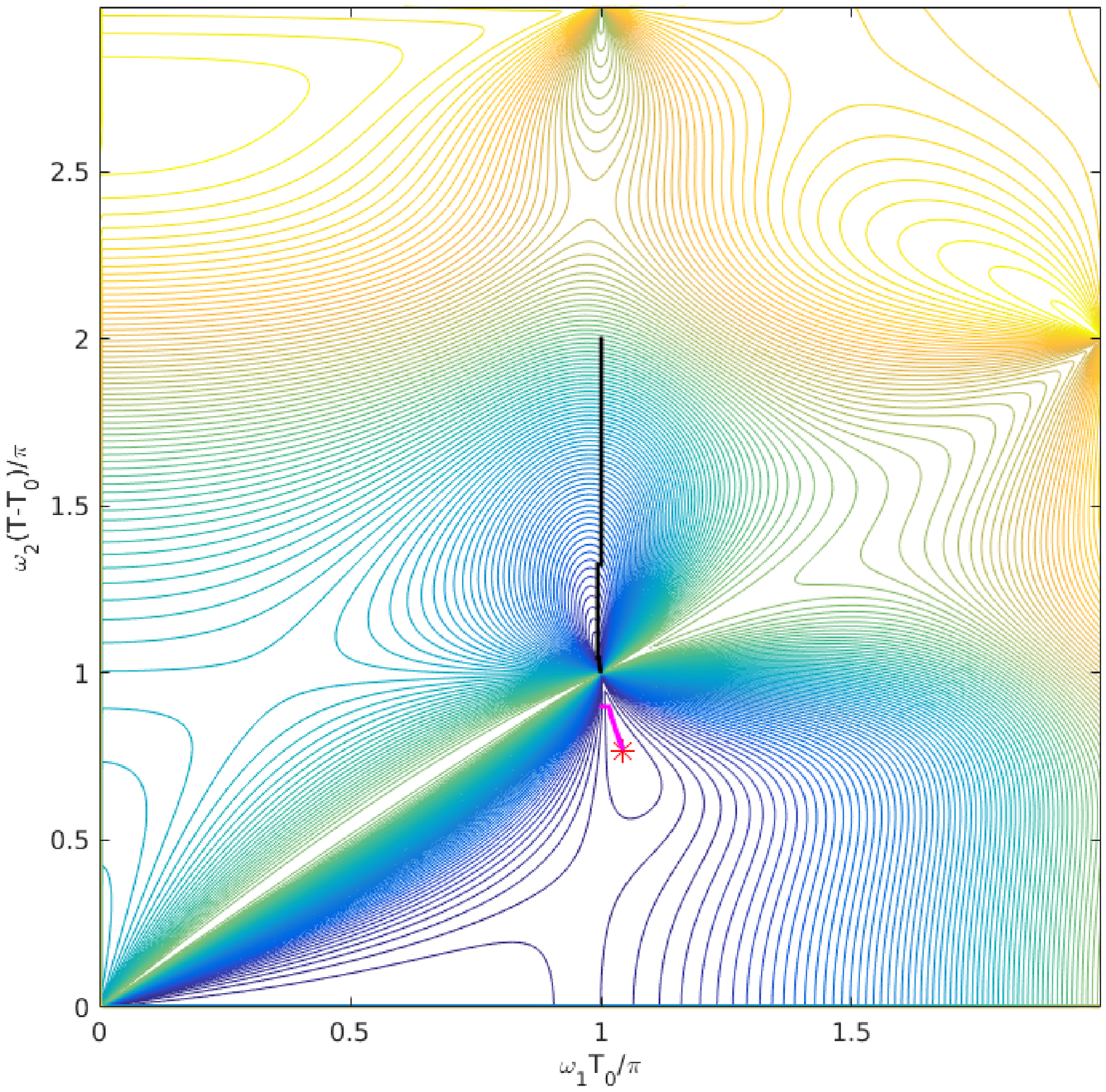}

\caption{\textbf{Up-Left:} A dog aortic pressure cycle (in blue), its dicrotic
notch (black dot), and the IMF (in red). \textbf{Up-Right:} heat-map
plot of $P\left(\frac{\omega_{1}T_{0}}{\pi},\frac{\omega_{2}\left(T-T_{0}\right)}{\pi}\right)$
for the cycle in left with the location of the solution marked with
red dot. \textbf{Down-Left:} Contour plot of $P\left(\frac{\omega_{1}T_{0}}{\pi},\frac{\omega_{2}\left(T-T_{0}\right)}{\pi}\right)$.
The location of the minimizer of $P$ is shown by a red dot. \textbf{Down-Right:}
Contour plot of $P\left(\frac{\omega_{1}T_{0}}{\pi},\frac{\omega_{2}\left(T-T_{0}\right)}{\pi}\right)$
and the location of the minimizer of $P$ tracked by the pattern search.
The true optimum point is marked with a red asterisk. The lower pattern
search set (in black) has converged towards the correct optimum. The
upper pattern search set (in magenta) has converged to a local minima
near the node. }

\label{Fig: Fit plot and Mesh 2}
\end{figure}

Before moving on, we show the topology of the $P\left(\omega_{1},\omega_{2}\right)$
function and also its minima locations in $\omega_{1}$ and $\omega_{2}$
space. These will provide useful insights on where to set the initialization
point(s) of a possible fast IF algorithm. The data description is
provided in the next section. In Figures \ref{Fig: Fit plot and Mesh 1}
and \ref{Fig: Fit plot and Mesh 2}, we have presented two different
dog aortic pressure cycles with the IMF extracted by the means of
the brute-force IF Algorithm \ref{Alg: InF Algorithm}. Figures \ref{Fig: Fit plot and Mesh 1}
and \ref{Fig: Fit plot and Mesh 2}, top right, show the heat-map
plots of $P\left(\frac{\omega_{1}T_{0}}{\pi},\frac{\omega_{2}\left(T-T_{0}\right)}{\pi}\right)$.
The complex nature of $P\left(\omega_{1},\omega_{2}\right)$ can be
seen in these figures. We purposefully plotted $P$ in the dimensionless
coordinates $\frac{\omega_{1}T_{0}}{\pi}$ and $\frac{\omega_{2}\left(T-T_{0}\right)}{\pi}$
to show the behavior of this function with respect to the lattice
node locations $\mathcal{N}$ defined in (\ref{eq: Nodes})-(\ref{eq: Lattice Nodes2}).
To have a better view and understanding of the $P\left(\omega_{1},\omega_{2}\right)$
topology, a contour of $P\left(\frac{\omega_{1}T_{0}}{\pi},\frac{\omega_{2}\left(T-T_{0}\right)}{\pi}\right)$
is shown in those figures. The general topology of $P\left(\frac{\omega_{1}T_{0}}{\pi},\frac{\omega_{2}\left(T-T_{0}\right)}{\pi}\right)$,
for all aortic or carotid pulse waveforms, is similar to the ones
presented in Figures \ref{Fig: Fit plot and Mesh 1} and \ref{Fig: Fit plot and Mesh 2}.
However, the location of the minimizer is not similar. 

Our investigations show that the locations of the minimizers of all
$P$ functions construct two different areas in the dimensionless
coordinates $\frac{\omega_{1}T_{0}}{\pi}$ and $\frac{\omega_{2}\left(T-T_{0}\right)}{\pi}$.
We call these areas as the \emph{upper lobe} and \emph{lower lobe}.
The upper lobe is an area, in $\mathcal{D}$, confined above the line
$\frac{\omega_{2}\left(T-T_{0}\right)}{\pi}=1$. The lower lobe is
an area, in $\mathcal{D}$, confined below the line $\frac{\omega_{2}\left(T-T_{0}\right)}{\pi}=1$.
This is also the case for human subject data \cite{pahlevan2017noninvasive}.
This type of topology suggests two critical initial guess areas for
any non-brute-force algorithm solving (\ref{eq: Discrete min}): one
set of points in the upper lobe, the other in the lower. In the remaining
part of this section, we introduce a fast IF algorithm based on the
pattern search method \cite{kolda2003optimization}. 

\subsection{Pattern Search IF}

The pattern search algorithm (or sometimes called the \emph{compass
search} algorithm) is explained in detail in \cite{kolda2003optimization}.
For completeness, we have summarized the pattern search algorithm
in Algorithm \ref{Alg: Pattern Search}. The convergence analysis
of this method is expressed in \cite{kolda2003optimization}.

\begin{algorithm}[H]
\begin{raggedright}
\textbf{Initialization.}
\par\end{raggedright}
\begin{raggedright}
Let $f:\mathbb{R}^{n}\rightarrow\mathbb{R}$ be given. 
\par\end{raggedright}
\begin{raggedright}
Let $x_{0}\in\mathbb{R}^{n}$ be the initial guess. 
\par\end{raggedright}
\begin{raggedright}
Let $\triangle_{tol}>0$ be the tolerance used to test for convergence. 
\par\end{raggedright}
\begin{raggedright}
Let $\triangle_{0}>\triangle_{tol}$ be the initial value of the step
length control parameter. 
\par\end{raggedright}
\begin{raggedright}
\textbf{Algorithm.} For each iteration $k=1,2,\ldots$
\par\end{raggedright}
\begin{raggedright}
\textbf{Step 1. }Let $\mathcal{D}_{\oplus}$ be the set of coordinate
directions $\left\{ \pm e_{i}\left|i=1,\ldots,n\right.\right\} $,
where $e_{i}$ is the \emph{i}th unit coordinate vector in $\mathbb{R}^{n}$. 
\par\end{raggedright}
\begin{raggedright}
\textbf{Step 2.} If there exists $d_{k}\in\mathcal{D}_{\oplus}$ such
that $f\left(x_{k}+\triangle_{k}d_{k}\right)<f\left(x_{k}\right)$,
then do the following: 
\par\end{raggedright}
\begin{itemize}
\item \begin{raggedright}
Set $x_{k+1}=x_{k}+\triangle_{k}d_{k}$.
\par\end{raggedright}
\item \begin{raggedright}
Set $\triangle_{k+1}=\triangle_{k}$.
\par\end{raggedright}
\end{itemize}
\begin{raggedright}
\textbf{Step 3.} Otherwise, $f\left(x_{k}+\triangle_{k}d_{k}\right)\geqslant f\left(x_{k}\right)$
for all $d_{k}\in\mathcal{D}_{\oplus}$, so do the following:
\par\end{raggedright}
\begin{itemize}
\item \begin{raggedright}
Set $x_{k+1}=x_{k}$.
\par\end{raggedright}
\item \begin{raggedright}
Set $\triangle_{k+1}=\frac{1}{2}\triangle_{k}$.
\par\end{raggedright}
\item \begin{raggedright}
If $\triangle_{k+1}<\triangle_{tol}$, then \textbf{terminate}.
\par\end{raggedright}
\end{itemize}
\caption{Pattern Search \cite{kolda2003optimization}}

\label{Alg: Pattern Search}
\end{algorithm}

The fast IF algorithm, without considering the nodes (\ref{eq: Nodes}),
is expressed in Algorithm \ref{Alg: Pattern Search IF}. As mentioned
before, what makes Algorithm \ref{Alg: Pattern Search IF} fast is
embedded in three different objects:
\begin{enumerate}
\item The initial guess set up in the initialization part of the algorithm. 
\item The fast analytic solution at each point iteration defined by (\ref{eq: P definition})
and (\ref{eq: Inner loop argmin sol}).
\item The pattern search part which is a substitute for the brute force
algorithm.
\end{enumerate}
Figure \ref{Fig: Fit plot and Mesh 1}, bottom right, shows the results
of Algorithm \ref{Alg: Pattern Search IF}. In this figure, when using
Algorithm \ref{Alg: Pattern Search IF}, we have used two initial
guesses $\left(\frac{\omega_{1}T_{0}}{\pi}=1,\frac{\omega_{2}\left(T-T_{0}\right)}{\pi}=2\right)$
and $\left(\frac{\omega_{1}T_{0}}{\pi}=1,\frac{\omega_{2}\left(T-T_{0}\right)}{\pi}=0.9\right)$.
As depicted on the figure, the initial guess located in the upper
lobe has converged towards the true minimizer in $\mathcal{D}$. On
a PC having 8 threads, Intel\textregistered{} Core\texttrademark{}
i7-4700MQ CPU @ 2.40GHz \texttimes{} 8, running a Matlab implementation
of the brute-force Algorithm \ref{Alg: InF Algorithm} in parallel
takes roughly $85$ seconds. On the other hand, achieving the same
minimizer, using a sequential version of the fast Algorithm \ref{Alg: Pattern Search IF},
takes approximately $0.5$ seconds. 

The same test was done for another aortic cycle presented in Figure
\ref{Fig: Fit plot and Mesh 2}. We used the same initial guesses
as before. This time, on the same PC, using the same implementations,
the brute-force Algorithm \ref{Alg: InF Algorithm} took roughly $80$
seconds and the fast Algorithm \ref{Alg: Pattern Search IF} took
approximately $0.5$ seconds. These two examples show a speed up of
almost $160$ times. In the next section we present more about the
statistical accuracy of Algorithm \ref{Alg: Pattern Search IF} and
its physiological capabilities. 

\begin{algorithm}[H]
\begin{raggedright}
\textbf{Initialization.}
\par\end{raggedright}
\begin{raggedright}
Let $\mathbf{f}\in\mathbb{R}^{n+m}$ be a given discrete aortic/carotid
signal with specified $T_{0}$ and $T-T_{0}$.
\par\end{raggedright}
\begin{raggedright}
Let $\mathcal{D}=\left\{ \left(\omega_{1},\omega_{2}\right)\left|0.5\leqslant\frac{\omega_{1}T_{0}}{\pi}\leqslant1.5,\,0.5\leqslant\frac{\omega_{2}\left(T-T_{0}\right)}{\pi}\leqslant3\right.\right\} $.
\par\end{raggedright}
\begin{raggedright}
Let $G=\bigcup_{l=1}^{M}\left\{ \left(\omega_{1},\omega_{2}\right)_{l}\in\mathcal{D}\right\} $
be the set of $M$ random initial guesses excluding the nodes (\ref{eq: Nodes}). 
\par\end{raggedright}
\begin{raggedright}
Let $\triangle\omega_{tol}>0$ be the convergence tolerance. 
\par\end{raggedright}
\begin{raggedright}
Let $\triangle\omega_{0}>\triangle\omega_{tol}$ be the initial step
length. 
\par\end{raggedright}
\begin{raggedright}
Let $\tilde{\omega}^{k}=\left(\omega_{1}^{k},\omega_{2}^{k}\right)$
and 
\[
P\left(\tilde{\omega}^{k}\right)=\left\Vert \mathbf{Q}\left(\tilde{\omega}^{k},b_{1}^{*}\left(\tilde{\omega}^{k}\right),b_{2}^{*}\left(\tilde{\omega}^{k}\right);\mathbf{t}\right)+\bar{p}^{*}\left(\tilde{\omega}^{k}\right)\mathbf{1}-\mathbf{f}\right\Vert _{2}^{2}
\]
 for the \emph{k}th iteration, defined by (\ref{eq: P definition}),
which is solved using (\ref{eq: Inner loop argmin sol}).
\par\end{raggedright}
\begin{raggedright}
Let $\mathcal{D}_{\oplus}=\left\{ \pm e_{j}\left|j=1,2\right.\right\} $,
where $e_{j}$ is the \emph{j}th unit coordinate vector in $\mathbb{R}^{2}$. 
\par\end{raggedright}
\begin{raggedright}
\textbf{Algorithm.} For each initial guess $\tilde{\omega}^{i}\in G$,
$i=1,\ldots,M$, and for each iteration $k_{i}=0,1,\ldots$ 
\par\end{raggedright}
\begin{raggedright}
\textbf{Step 1.} If there exists $\mathbf{d}_{k_{i}}\in\mathcal{D}_{\oplus}$
such that $P\left(\tilde{\omega}^{k_{i}}+\triangle\omega_{k_{i}}\mathbf{d}_{k_{i}}\right)<P\left(\tilde{\omega}^{k_{i}}\right)$,
then: 
\par\end{raggedright}
\begin{itemize}
\item \begin{raggedright}
$\tilde{\omega}^{k_{i}+1}=\tilde{\omega}^{k_{i}}+\triangle\omega_{k_{i}}\mathbf{d}_{k_{i}}$.
\par\end{raggedright}
\item \begin{raggedright}
$\triangle\omega_{k_{i}+1}=\triangle\omega_{k_{i}}$.
\par\end{raggedright}
\end{itemize}
\begin{raggedright}
\textbf{Step 2.} Otherwise, if $P\left(\tilde{\omega}^{k_{i}}+\triangle\omega_{k_{i}}\mathbf{d}_{k_{i}}\right)\geqslant P\left(\tilde{\omega}^{k_{i}}\right)$
for all $\mathbf{d}_{k_{i}}\in\mathcal{D}_{\oplus}$, then:
\par\end{raggedright}
\begin{itemize}
\item \begin{raggedright}
$\tilde{\omega}^{k_{i}+1}=\tilde{\omega}^{k_{i}}$.
\par\end{raggedright}
\item \begin{raggedright}
$\triangle\omega_{k_{i}+1}=\frac{1}{2}\triangle\omega_{k_{i}}$.
\par\end{raggedright}
\item \begin{raggedright}
If $\triangle\omega_{k_{i}+1}<\triangle\omega_{tol}$, then \textbf{terminate}
and $\tilde{\omega}_{i}^{*}=\tilde{\omega}^{k_{i}+1}$.
\par\end{raggedright}
\end{itemize}
\begin{raggedright}
\textbf{Step 3.} The solution is$\tilde{\omega}^{*}=\underset{i\in\left\{ 1,\ldots,M\right\} }{\arg\min}\,P\left(\tilde{\omega}_{i}^{*}\right)$. 
\par\end{raggedright}
\caption{Fast IF}

\label{Alg: Pattern Search IF}
\end{algorithm}

\section{Real Data Example}

The real dog data used in this manuscript is well described in \cite{swamy2009continuous}.
Since, at the time of the the data retrieval, the data was downloaded
with different sampling rates, we re-sampled all six dog data at $500\,Hz$.
We used a modified version of the automatic cycle selection introduced
in \cite{zong2003open} to pick cycles. Dicrotic notch locations were
then found from the picked cycles \cite{li2010automatic}. We totally
extracted $59384$ acceptable aortic cycles form those six dogs. 

\subsection{Statistical Accuracy}

To check the statistical accuracy of the fast IF algorithm versus
the brute-force IF algorithm, we compared the results of these two
algorithms on the extracted $59384$ dog aortic cycles. The brute-force
IF algorithm (Algorithm \ref{Alg: InF Algorithm}) was run over the
sample set with a mesh size ($\underset{l\neq m}{\min}\left(\omega_{1}^{l}-\omega_{1}^{m}\right)=\underset{l\neq m}{\min}\left(\omega_{2}^{l}-\omega_{2}^{m}\right)$)
of $0.02\pi$. Algorithm \ref{Alg: Pattern Search IF} was run on
the same sample set of $59384$ aortic cycles with $\triangle\omega_{tol}=0.001$,
and $\triangle\omega_{0}=0.1$, comprising a mesh size of $\frac{0.1}{2^{6}}$.
The brute-force algorithm has a larger mesh size due to heavy computational
cost of this algorithm. The maximum average difference between the
IFs found by these two algorithms was found to be less than $0.0475$.
This difference is smaller than both mesh sizes used for the brute-force
and fast IF algorithms. This shows that, on average, the fast IF algorithm
(Algorithm \ref{Alg: Pattern Search IF}) reaches the same minima
as the brute-force algorithm (Algorithm \ref{Alg: InF Algorithm}). 

\subsection{Physiological Observations}

\begin{figure}[H]
$ $\includegraphics[scale=0.5]{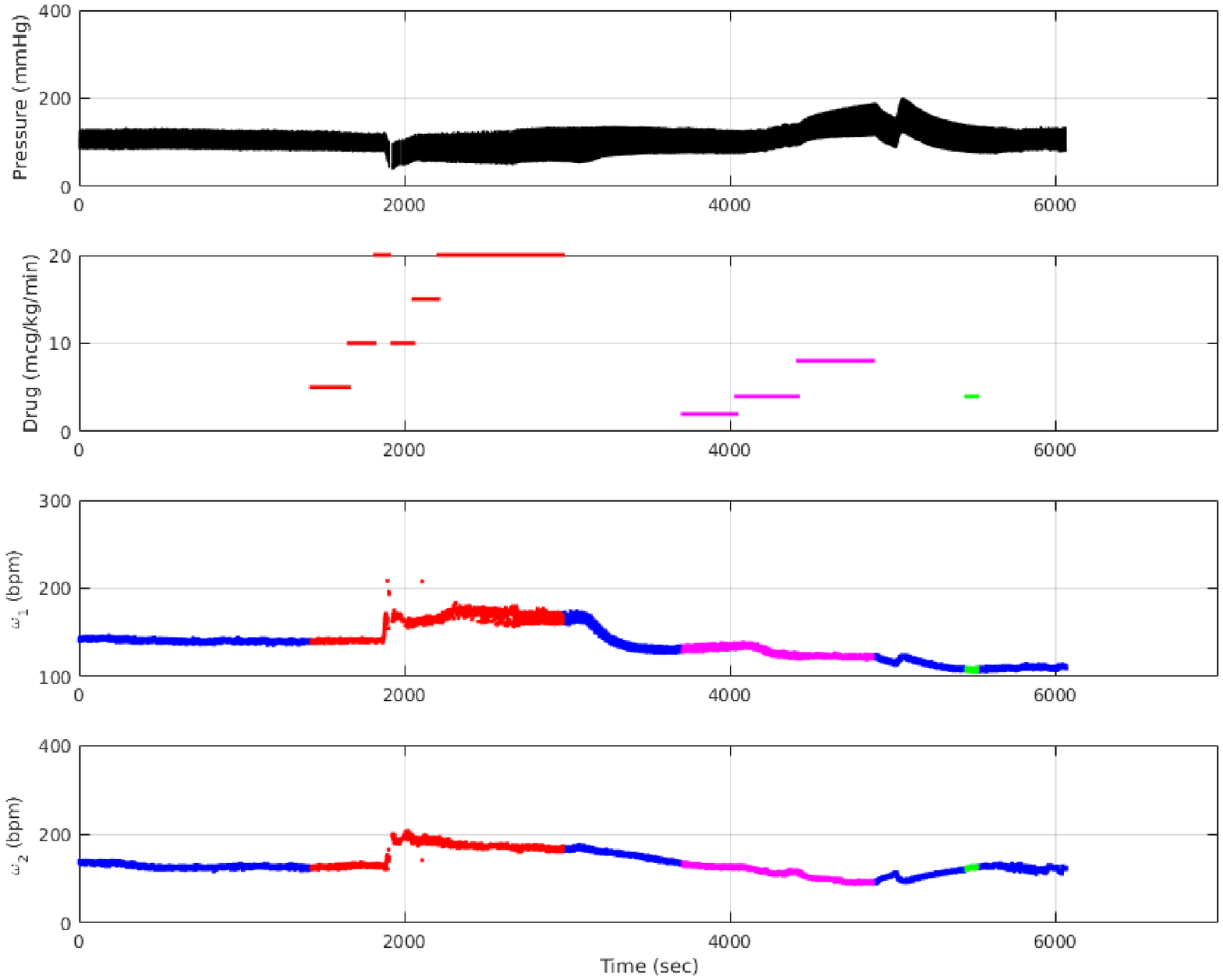}

\caption{Drug effects on $\omega_{1}$ and $\omega_{2}$. \textbf{First Panel:}
The measured aortic pressure waveform recorded in time. \textbf{Second
Panel:} Dosage of dobutamine (in red), phenylephrine (in purple),
and nitroglycerin (in green) during the aortic pressure measurement.
\textbf{Third Panel:} Changes of $\omega_{1}$ in units of bit per
minute (bpm) over the measurement time. Each drug effect is projected
with its corresponding color. No drug areas are in blue. \textbf{Fourth
Panel:} Changes of $\omega_{2}$ in units of bpm over the measurement
time. Each drug effect is projected with its corresponding color.
No drug areas are in blue.}

\label{Fig: Drugs}
\end{figure}

To evaluate the new fast IF algorithm (Algorithm \ref{Alg: Pattern Search IF}),
we applied the algorithm on the measured aortic pressure signal from
one dog experiencing various pharmacological interventions, see Figure
\ref{Fig: Drugs}. During the experiment, the dog was under the following
pharmacological influences: infusion of dobutamine (5-20 $\mu g/kg/min$),
phenylephrine (2-8 $\mu g/kg/min$) and nitroglycerin (4 $\mu g/kg/min$)
during different time intervals. 

The third panel, in Figure \ref{Fig: Drugs}, shows the dosage and
duration of each drug in the experiment. In the first phase of the
experiment dobutamine has been injected at a low dosage followed by
a fluctuation in the dosage of injection. The effect of dobutamine
on the cardiovascular system is to increase the strength and force
of the heartbeat. Consequently, it forces more blood to circulate
throughout the body. In previous works \cite{pahlevan2017noninvasive,pahlevan2014intrinsic},
we hypothesized that $\omega_{1}$ would be a representative of heart
functionality. We also hypothesized that $\omega_{1}$ and $\omega_{2}$
would try to keep a balance during changes. These hypotheses can be
seen during the injection of dobutamine in this figure.

Next, phenylephrine has been injected at a low dosage and the dosage
is then increased over time. Phenylephrine is a decongestant which
affects the cardiovascular system by shrinking blood vessels. $\omega_{2}$
shows an almost monotone decrease during the infusion of phenylephrine.
This is again in qualitative accord with what we presented in \cite{pahlevan2014intrinsic}.

Lastly, nitroglycerin has been injected at a constant dosage. Nitroglycerin
helps to dilate the blood vessels. This dilation can be captured with
$\omega_{2}$, as can be seen from the figure. Generally, based on
this figure, IFs are able to capture changes in the dynamics of the
system under the effects of different drugs.

\section{Conclusion}

In this paper, we provided a mathematical foundation for the IF model
\cite{tavallali2015convergence}. We showed how to derive an estimation
of the IF model (\ref{eq: The Trend}) by considering basic physics
principles. More precisely, we showed that the IF model can be estimated
from Navier-Stokes and elasticity equations. 

We further analysed the IF model (\ref{eq: Discrete min}). This helped
to introduce a fast algorithm for the IF method (Algorithm \ref{Alg: Pattern Search IF}).
What made this algorithm fast was embedded in the proper set up of
the initial guesses based on the topology of the problem, fast analytic
solution at each point iteration, and substituting the brute force
algorithm with a pattern search method. These changes would convert
an iterative and brute-force method (Algorithm \ref{Alg: InF Algorithm})
into an algebraic and iterative method (Algorithm \ref{Alg: Pattern Search IF}).
The presented fast algorithm, in this article, has a speed up of more
than 100 times compared to the brute-force algorithm provided in \cite{tavallali2015convergence}.
From a statistical perspective, we have also shown that the algorithm
presented in this article complies well with the brute-force implementations
of this method.

We also showed, on a real dataset, that the fast IF Algorithm \ref{Alg: Pattern Search IF}
can depict correlations between its outputs and infusion of certain
drugs. This part of our paper can be subject to further physiological
and clinical investigations in a future work.

\section{Authors' Contributions}

\textcolor{black}{P.T. conceived of the mathematical and numerical
methods of the study, carried out the modeling, programmed the initial
code of the method, and drafted the manuscript; H.K. helped with the
mathematical and numerical derivations, helped draft and revise the
manuscript, and conducted the real data case example; J.K. conducted
the brute-force simulations on the real data case example, and helped
draft the manuscript. All authors gave final approval for publication.}

\section{Acknowledgement}

We would like to thank Mr. Sean Brady for constructive discussions
and editorial comments. 

\section{Research Ethics}

All experiments and procedures were reviewed and approved by the MSU
All-University Committee on Animal Use and Care \cite{swamy2009continuous}.

\section{Permission to Carry Out Fieldwork}

This study did not have fieldwork.

\section{Funding}

This work was not funded.

\section{Figures}

\section*{Appendix A\label{sec:Appendix-A}}

In this appendix, we show how one can derive (\ref{eq: Simplified Elasticity})
and (\ref{eq: Linearized Momentum}), from the Navier-Stokes and elasticity
equations. \textcolor{black}{Having $x$, $r$ and $\theta$ as the
cylindrical coordinate system, with $x$ in the direction of the aortic
length, the momentum and the continuity equations are} 
\begin{equation}
\begin{array}{c}
\rho\left(\frac{\partial u}{\partial t}+u\frac{\partial u}{\partial x}+v\frac{\partial u}{\partial r}+\frac{\omega}{r}\frac{\partial u}{\partial\theta}\right)+\frac{\partial P}{\partial x}=\\
\mu\left(\frac{\partial^{2}u}{\partial x^{2}}+\frac{\partial^{2}u}{\partial r^{2}}+\frac{1}{r}\frac{\partial u}{\partial r}+\frac{1}{r^{2}}\frac{\partial^{2}u}{\partial\theta^{2}}\right),
\end{array}\label{eq:7}
\end{equation}
\begin{equation}
\begin{array}{c}
\rho\left(\frac{\partial v}{\partial t}+u\frac{\partial v}{\partial x}+v\frac{\partial v}{\partial r}+\frac{w}{r}\frac{\partial v}{\partial\theta}-\frac{w^{2}}{r}\right)+\frac{\partial P}{\partial r}=\\
\mu\left(\frac{\partial^{2}v}{\partial x^{2}}+\frac{\partial^{2}v}{\partial r^{2}}+\frac{1}{r}\frac{\partial v}{\partial r}-\frac{v}{r^{2}}+\frac{1}{r^{2}}\frac{\partial^{2}v}{\partial\theta^{2}}-\frac{2}{r^{2}}\frac{\partial w}{\partial\theta}\right),
\end{array}\label{eq:8}
\end{equation}
\begin{equation}
\begin{array}{c}
\rho\left(\frac{\partial w}{\partial t}+u\frac{\partial w}{\partial x}+v\frac{\partial w}{\partial r}+\frac{w}{r}\frac{\partial w}{\partial\theta}+\frac{vw}{r}\right)+\frac{1}{r}\frac{\partial P}{\partial\theta}=\\
\mu\left(\frac{\partial^{2}w}{\partial x^{2}}+\frac{\partial^{2}w}{\partial r^{2}}+\frac{1}{r}\frac{\partial w}{\partial r}-\frac{w}{r^{2}}+\frac{1}{r^{2}}\frac{\partial^{2}w}{\partial\theta^{2}}+\frac{2}{r^{2}}\frac{\partial v}{\partial\theta}\right),
\end{array}\label{eq:9}
\end{equation}
\begin{equation}
\frac{\partial u}{\partial x}+\frac{\partial v}{\partial r}+\frac{v}{r}+\frac{1}{r}\frac{\partial w}{\partial\theta}=0.\label{eq:10}
\end{equation}
Here, $u$, $v$ and $w$ are velocity vector components in $x$,
$r$ and $\theta$ directions, respectively. We assume that aorta
is a straight and sufficiently long tube with constant circular cross
section with the tube wall following classical elasticity theory dynamics
and blood is considered to be an incompressible Newtonian fluid with
the velocity field being axisymmetric. In the absents of any external
forces causing flow rotation, the assumption that the flow field is
symmetrical about the longitudinal axis of the tube is justified.
This means $w=\frac{\partial w}{\partial\theta}=\frac{\partial v}{\partial\theta}=\frac{\partial u}{\partial\theta}=\frac{\partial P}{\partial\theta}=0$.
Hence, Equations (\ref{eq:7})-(\ref{eq:10}) will be simplified as
\begin{equation}
\rho\left(\frac{\partial u}{\partial t}+u\frac{\partial u}{\partial x}+v\frac{\partial u}{\partial r}\right)+\frac{\partial P}{\partial x}=\mu\left(\frac{\partial^{2}u}{\partial x^{2}}+\frac{\partial^{2}u}{\partial r^{2}}+\frac{1}{r}\frac{\partial u}{\partial r}\right),\label{eq:11}
\end{equation}
\begin{equation}
\rho\left(\frac{\partial v}{\partial t}+u\frac{\partial v}{\partial x}+v\frac{\partial v}{\partial r}\right)+\frac{\partial P}{\partial r}=\mu\left(\frac{\partial^{2}v}{\partial x^{2}}+\frac{\partial^{2}v}{\partial r^{2}}+\frac{1}{r}\frac{\partial v}{\partial r}-\frac{v}{r^{2}}\right),\label{eq:12}
\end{equation}
\begin{equation}
\frac{\partial u}{\partial x}+\frac{\partial v}{\partial r}+\frac{v}{r}=0.\label{eq:13}
\end{equation}
Since the radius of the tube $a$ is smaller than the length of the
tube, $a\ll h$, and also the average velocity of the blood in aorta
is smaller than the speed of wave propagation \cite{zamir2000physics},
we have
\begin{equation}
u\frac{\partial u}{\partial x}\ll\frac{\partial u}{\partial t},\label{eq:14}
\end{equation}
\begin{equation}
v\frac{\partial u}{\partial r}\ll\frac{\partial u}{\partial t},\label{eq:15}
\end{equation}
\begin{equation}
u\frac{\partial v}{\partial x}\ll\frac{\partial v}{\partial t},\label{eq:16}
\end{equation}
\begin{equation}
v\frac{\partial v}{\partial r}\ll\frac{\partial v}{\partial t},\label{eq:17}
\end{equation}
\begin{equation}
\frac{\partial^{2}u}{\partial x^{2}}\ll\frac{\partial^{2}u}{\partial r^{2}},\label{eq:18}
\end{equation}
\begin{equation}
\frac{\partial^{2}v}{\partial x^{2}}\ll\frac{\partial^{2}v}{\partial r^{2}}.\label{eq:19}
\end{equation}
Using these, Equations (\ref{eq:11})-(\ref{eq:13}) will reduce to
\begin{equation}
\rho\frac{\partial u}{\partial t}+\frac{\partial P}{\partial x}=\mu\left(\frac{\partial^{2}u}{\partial r^{2}}+\frac{1}{r}\frac{\partial u}{\partial r}\right)\label{eq:20}
\end{equation}
\begin{equation}
\rho\frac{\partial v}{\partial t}+\frac{\partial P}{\partial r}=\mu\left(\frac{\partial^{2}v}{\partial r^{2}}+\frac{1}{r}\frac{\partial v}{\partial r}-\frac{v}{r^{2}}\right)\label{eq:21}
\end{equation}
\begin{equation}
\frac{\partial u}{\partial x}+\frac{\partial v}{\partial r}+\frac{v}{r}=0\label{eq:22}
\end{equation}
It is important to mention that the velocity vector is a function
of time and location. In other words, we have $u\left(x,r,t\right)$
and $v\left(x,r,t\right)$. Based on the characteristic length of
the problem, pressure $P$ can be assumed to be a function of $x$
and $t$ and not $r$, i.e. $P\left(x,t\right)$ (See Chapter 5 of
\cite{zamir2000physics}). Hence we can set $\frac{\partial P}{\partial r}\approx0$
in (\ref{eq:21}).

\subsection*{Momentum Equations:}

Considering (\ref{eq:20}), one can integrate both sides with respect
to the differential element of the area $2\pi rdr$. 
\begin{equation}
\intop_{0}^{a\left(x,t\right)}2\pi r\rho\frac{\partial u}{\partial t}dr+\intop_{0}^{a\left(x,t\right)}2\pi r\frac{\partial P}{\partial x}dr=\intop_{0}^{a\left(x,t\right)}2\pi r\mu\left(\frac{\partial^{2}u}{\partial r^{2}}+\frac{1}{r}\frac{\partial u}{\partial r}\right)dr.\label{eq:23}
\end{equation}
The upper boundary of this integral is the radius $a\left(x,t\right)$
of the tube. We know that the flow $Q\left(x,t\right)$ is defined
as 
\begin{equation}
Q\left(x,t\right)=\intop_{0}^{a\left(x,t\right)}2\pi ru\left(x,r,t\right)dr.\label{eq:24}
\end{equation}
Hence, using Leibniz rule we can find the derivative of the flow with
respect to time as 
\begin{equation}
\frac{\partial Q}{\partial t}\left(x,t\right)=\intop_{0}^{a\left(x,t\right)}2\pi r\frac{\partial u}{\partial t}\left(x,r,t\right)dr+2\pi a\left(x,t\right)\frac{\partial a}{\partial t}\left(x,t\right)u\left(x,a\left(x,t\right),t\right).\label{eq:25}
\end{equation}
Considering no slip boundary condition $u\left(x,a\left(x,t\right),t\right)=0$
on the tube wall, Equation (\ref{eq:25}) reduces to 
\begin{equation}
\frac{\partial Q}{\partial t}\left(x,t\right)=\intop_{0}^{a\left(x,t\right)}2\pi r\frac{\partial u}{\partial t}\left(x,r,t\right)dr.\label{eq:26}
\end{equation}
Therefore Equation (\ref{eq:23}) will be simplified to
\begin{equation}
\rho\frac{\partial Q}{\partial t}+\pi a^{2}\frac{\partial P}{\partial x}=2\pi\mu\intop_{0}^{a\left(x,t\right)}r\left(\frac{\partial^{2}u}{\partial r^{2}}+\frac{1}{r}\frac{\partial u}{\partial r}\right)dr.\label{eq:27}
\end{equation}
Since we have 
\begin{equation}
\intop_{0}^{a\left(x,t\right)}r\left(\frac{\partial^{2}u}{\partial r^{2}}+\frac{1}{r}\frac{\partial u}{\partial r}\right)dr=\intop_{0}^{a\left(x,t\right)}\frac{\partial\left(r\frac{\partial u}{\partial r}\right)}{\partial r}dr,
\end{equation}
Equation (\ref{eq:27}) would become

\begin{equation}
\rho\frac{\partial Q}{\partial t}+\pi a^{2}\frac{\partial P}{\partial x}=2\pi\mu r\frac{\partial u}{\partial r}\mid_{0}^{a\left(x,t\right)}.\label{eq:55}
\end{equation}
Using separation of variables $u\left(x,r,t\right)=U\left(r\right)\bar{u}\left(x,t\right)$,
for some function $U\left(r\right)$ and $\bar{u}\left(x,t\right)$,
Equation (\ref{eq:55}) will be simplified to 

\begin{equation}
\rho\frac{\partial Q}{\partial t}+\pi a^{2}\frac{\partial P}{\partial x}=2\pi\mu a\left(x,t\right)\bar{u}\left(x,t\right)\left(\frac{dU}{dr}\mid_{a\left(x,t\right)}\right).\label{eq:56}
\end{equation}
\textcolor{black}{From Equation (\ref{eq:24}) and }$u\left(x,r,t\right)=U\left(r\right)\bar{u}\left(x,t\right)$,
we have

\begin{equation}
\bar{u}\left(x,t\right)=Q\left(x,t\right)\left(2\pi\intop_{0}^{a\left(x,t\right)}rU\left(r\right)dr\right)^{-1}.\label{eq:58}
\end{equation}
Therefore equation (\ref{eq:56}) will be simplified to 

\begin{equation}
\rho\frac{\partial Q}{\partial t}\left(x,t\right)+\pi a^{2}\left(x,t\right)\frac{\partial P}{\partial x}\left(x,t\right)=\mu a\left(x,t\right)\left(\frac{dU}{dr}\mid_{a\left(x,t\right)}\right)\left(\intop_{0}^{a\left(x,t\right)}rU\left(r\right)dr\right)^{-1}Q\left(x,t\right).\label{eq:59}
\end{equation}
From (\ref{eq:59}), we can relabel some terms and introduce the inductance
$\mathcal{L}$ and resistance $\mathcal{R}$ as

\begin{equation}
\mathcal{L}\left(x,t\right)=\frac{\rho}{\pi a^{2}\left(x,t\right)},\label{eq:60}
\end{equation}

\begin{equation}
\mathcal{R}\left(x,t\right)=-\mu\left(\frac{dU}{dr}\mid_{a\left(x,t\right)}\right)\left(\pi a\left(x,t\right)\intop_{0}^{a\left(x,t\right)}rU\left(r\right)dr\right)^{-1}.\label{eq:61}
\end{equation}
These will convert (\ref{eq:59}) into 
\begin{equation}
-\frac{\partial P}{\partial x}\left(x,t\right)=\mathcal{L}\left(x,t\right)\frac{\partial Q}{\partial t}\left(x,t\right)+\mathcal{R}\left(x,t\right)Q\left(x,t\right).\label{eq: Linearized Momentum nonConst}
\end{equation}

\subsection*{Continuity Equation:}

Again, applying Leibniz rule to equation (\ref{eq:24}) we can find
the derivative of the flow with respect to $x$ 
\begin{equation}
\frac{\partial Q}{\partial x}\left(x,t\right)=\intop_{0}^{a\left(x,t\right)}2\pi r\frac{\partial u}{\partial x}\left(x,r,t\right)dr+2\pi a\left(x,t\right)\frac{\partial a}{\partial x}\left(x,t\right)u\left(x,a\left(x,t\right),t\right).\label{eq:33}
\end{equation}
Considering the no slip boundary condition $u\left(x,a(x,t),t\right)=0$,
equation (\ref{eq:33}) reduces to 
\begin{equation}
\frac{\partial Q}{\partial x}\left(x,t\right)=\intop_{0}^{a\left(x,t\right)}2\pi r\frac{\partial u}{\partial x}\left(x,r,t\right)dr.\label{eq:34}
\end{equation}
Now, we can rewrite the equation of continuity (\ref{eq:22}) as 
\begin{equation}
\intop_{0}^{a\left(x,t\right)}2\pi r\frac{\partial u}{\partial x}dr+\intop_{0}^{a\left(x,t\right)}2\pi\left(r\frac{\partial v}{\partial r}+v\right)dr=0.\label{eq:32}
\end{equation}
This equation, using (\ref{eq:34}), will result in
\begin{equation}
\frac{\partial Q}{\partial x}+\intop_{0}^{a\left(x,t\right)}2\pi\frac{\partial\left(rv\right)}{\partial r}dr=0.\label{eq:35}
\end{equation}
Simplifying the latter would show that 
\begin{equation}
\frac{\partial Q}{\partial x}+2\pi a\left(x,t\right)v\left(a\left(x,t\right)\right)=0.\label{eq:36}
\end{equation}
We note that $v\left(a\left(x,t\right)\right)=\frac{\partial a}{\partial t}\left(x,t\right)$.
Hence, having $A=\pi a^{2}\left(x,t\right)$, we can conclude $\frac{\partial A}{\partial t}=2\pi a\left(x,t\right)v\left(a\left(x,t\right)\right)$.
Consequently, Equation (\ref{eq:36}) can be written as
\begin{equation}
\frac{\partial Q}{\partial x}+\frac{\partial A}{\partial t}=0.\label{eq:37}
\end{equation}
Using the chain rule we have 
\begin{equation}
\frac{\partial A}{\partial t}=\frac{\partial A}{\partial P}\frac{\partial P}{\partial t}.\label{eq:Chain}
\end{equation}
Considering the wave speed $c_{0}$ of an incompressible fluid in
an elastic tube we have $\frac{\partial A}{\partial P}=\frac{A}{\rho c_{0}^{2}}$
\cite{parker2009brief}. Matching this with Equations (\ref{eq:37})
and (\ref{eq:35}) would result in
\begin{equation}
\frac{\partial Q}{\partial x}\left(x,t\right)+\frac{\pi a^{2}\left(x,t\right)}{\rho c_{0}^{2}}\frac{\partial P}{\partial t}\left(x,t\right)=0.\label{eq:38}
\end{equation}
In this equation, we can relabel $\frac{\pi a^{2}\left(x,t\right)}{\rho c_{0}^{2}}$
as the compliance $\mathcal{C}\left(x,t\right)$. Hence, (\ref{eq:38})
would become 
\begin{equation}
-\frac{\partial Q}{\partial x}\left(x,t\right)=\mathcal{C}\left(x,t\right)\frac{\partial P}{\partial t}\left(x,t\right).\label{eq: Continuity nonCnst}
\end{equation}

\subsection*{Wave Equations:}

As depicted so far, using the approximations in this appendix, and
considering the mentioned assumptions, we can characterize the wave
dynamics of the blood flow in aorta using the hyperbolic equations
(\ref{eq: Continuity nonCnst}) and (\ref{eq: Linearized Momentum nonConst}),
namely 
\begin{equation}
-\frac{\partial Q}{\partial x}\left(x,t\right)=\mathcal{C}\left(x,t\right)\frac{\partial P}{\partial t}\left(x,t\right),
\end{equation}
\begin{equation}
-\frac{\partial P}{\partial x}\left(x,t\right)=\mathcal{L}\left(x,t\right)\frac{\partial Q}{\partial t}\left(x,t\right)+\mathcal{R}\left(x,t\right)Q\left(x,t\right).
\end{equation}
The coefficients $\mathcal{C}\left(x,t\right)$, $\mathcal{L}\left(x,t\right)$
and $\mathcal{R}\left(x,t\right)$, in Equations (\ref{eq: Continuity nonCnst})
and (\ref{eq: Linearized Momentum nonConst}), are all positive and
functions of $a\left(x,t\right)$. However, since $a\left(x,t\right)$
is not changing drastically with respect to $x$ and $t$, we can
approximate all these coefficients with their corresponding constant
mean values $C$, $L$ and $R$. Using this approximation, we can
rewrite the wave equations as
\begin{equation}
-\frac{\partial Q}{\partial x}\left(x,t\right)=C\frac{\partial P}{\partial t}\left(x,t\right),
\end{equation}
\begin{equation}
-\frac{\partial P}{\partial x}\left(x,t\right)=L\frac{\partial Q}{\partial t}\left(x,t\right)+RQ\left(x,t\right).
\end{equation}

\end{document}